\documentclass[11pt]{amsart}
\usepackage{amssymb}
\usepackage{amsmath,amssymb}

\theoremstyle{plain}
\newtheorem{thm}{Theorem}[section]
\newtheorem{theorem}[thm]{Theorem}

\newtheorem{lemma}[thm]{Lemma}
\newtheorem{corollary}[thm]{Corollary}
\newtheorem{proposition}[thm]{Proposition}
\theoremstyle{definition}
\newtheorem{remark}[thm]{Remark}

\newtheorem{notation}[thm]{Notation}

\newtheorem{definition}[thm]{Definition}

\newtheorem{assumption}[thm]{Assumption}

\newtheorem{example}[thm]{Example}

\newtheorem{question}[thm]{Question}

\newtheorem{setup}[thm]{Setup}
\numberwithin{equation}{section}
\newcommand{\II}{{\rm II}}
\newcommand{\III}{{\rm III}}

\newcommand{\p}{\partial}

\newcommand{\sA}{{\mathcal A}}
\newcommand{\sB}{{\mathcal B}}
\newcommand{\sC}{{\mathcal C}}
\newcommand{\sD}{{\mathcal D}}
\newcommand{\sE}{{\mathcal E}}
\newcommand{\sF}{{\mathcal F}}

\newcommand{\sH}{{\mathcal H}}

\newcommand{\sK}{{\mathcal K}}

\newcommand{\sO}{{\mathcal O}}

\newcommand{\sS}{{\mathcal S}}
\newcommand{\sT}{{\mathcal T}}

\newcommand{\sV}{{\mathcal V}}

\newcommand{\C}{{\mathbb C}}

\newcommand{\N}{{\mathbb N}}
\newcommand{\BP}{{\mathbb P}}

\newcommand{\Q}{{\mathbb Q}}

\newcommand{\V}{{\mathbb V}}

\newcommand{\Z}{{\mathbb Z}}

\newcommand{\wg}{\widehat{\mathfrak l}}
\newcommand{\bV}{\mathbf{V}}

\newcommand{\bG}{\mathbf{G}}

\newcommand{\bK}{\mathbf{K}}
\newcommand{\bW}{\mathbf{W}}

\newcommand{\bg}{\mathbf{g}}

\newcommand{\fg}{{\mathfrak g}}
\newcommand{\fsl}{{\mathfrak s}{\mathfrak l}}
\newcommand{\fgl}{{\mathfrak g}{\mathfrak l}}
\newcommand{\ff}{{\mathfrak f}}

\newcommand{\fp}{{\mathfrak p}}
\newcommand{\fn}{{\mathfrak n}}

\newcommand{\fl}{{\mathfrak l}}

\newcommand{\fo}{{\mathfrak o}}
\newcommand{\fq}{{\mathfrak q}}

\newcommand{\fh}{{\mathfrak h}}

\newcommand{\fs}{{\mathfrak s}}
\newcommand{\ft}{{\mathfrak t}}

\newcommand{\aut}{{\mathfrak a}{\mathfrak u}{\mathfrak t}}

\newcommand\Aut{\rm Aut}

\newcommand{\bc}{\mathbf{c}}
\newcommand{\bw}{\mathbf{w}}

\def\Sym{\mathop{\rm Sym}\nolimits}

\def\Hom{\mathop{\rm Hom}\nolimits}

\title[Minimal rational curves and 1-flat G-structures]{Minimal rational curves and 1-flat irreducible  G-structures}

\author{Jun-Muk Hwang and Qifeng Li}

\thanks{This work was supported by the Institute for Basic Science (IBS-R032-D1).}

\begin{document}

\maketitle

\begin{center}
{\em Dedicated to the memory of Nessim Sibony. }
\end{center}

\begin{abstract}
1-flat irreducible G-structures, equivalently, irreducible G-structures admitting torsion-free affine connections, have been studied  extensively in differential geometry, especially in connection with the theory of affine holonomy groups. We propose to study them in a setting in algebraic geometry, where they arise from varieties of minimal rational tangents (VMRT) associated to families of minimal rational curves on  uniruled projective manifolds.
 We prove that such a structure is locally symmetric when the dimension of the uniruled projective manifold is at least 5. By the classification result of Merkulov and Schwachh\"ofer on irreducible affine holonomy, the problem is reduced to the case when the VMRT at a general point of the uniruled projective manifold is isomorphic to a subadjoint variety. In the latter situation, we prove a stronger result that, without the assumption of 1-flatness, the structure arising from VMRT is always locally flat.  The proof employs the method of  Cartan connections. An interesting feature is that Cartan connections are considered not for  the G-structures themselves, but for certain geometric structures on the spaces of minimal rational curves.
 \end{abstract}


\section{Introduction}\label{s.I}
We work in the holomorphic category: all varieties, groups and algebras are defined over complex numbers and all morphisms are holomorphic.  Unless mentioned otherwise, open sets refer to
Euclidean topology.  For a vector space $V$, its dual vector space is denoted by $V^{\vee}$. The projectivization $\BP
V$ of a vector space (or a vector bundle) $V$ is taken in the
classical sense, i.e., as the set of 1-dimensional subspaces of $V$. For a subvariety $Z \subset \BP V$,
its affine cone is denoted by $\widehat{Z} \subset V$. The group of linear automorphisms preserving $\widehat{Z}$ is denoted by $\Aut(\widehat{Z}) \subset {\rm GL}(V)$ and its Lie algebra is denoted by $\aut(\widehat{Z}) \subset \fgl(V)$.

\medskip
Recall that for a closed  subgroup $G \subset {\rm GL}(V)$ of the general linear group,
a $G$-structure on a manifold $Y$ with $\dim Y = \dim V$ is a $G$-principal subbundle of the frame bundle of $Y$. When we do not want to specify the structure group $G \subset {\rm GL}(V)$, we  just call it a G-structure (instead of $G$-structure with italic $G$).   G-structures are among basic objects of interest in differential geometry.
A fundamental problem is to find conditions for the local flatness of a G-structure, i.e., the condition for the existence of holomorphic charts on the manifold which send the G-structure to the flat G-structure on the vector space $V$ regarded as a manifold. There is also  the notion of flatness of finite order, a weaker version of local flatness. For example, a G-structure is 1-flat if it is flat up to the first jet at each point of the base manifold (see p. 528 of \cite{Br} for a precise definition).  The  1-flatness of a G-structure is equivalent to the existence of a torsion-free affine connection compatible with the G-structure (see p. 529 of \cite{Br}) and has been studied extensively in the theory of affine holonomy (see e.g. \cite{Br} and \cite{MS}). If the affine connection has a parallel curvature tensor, the 1-flat G-structure is said to be locally symmetric. Locally symmetric G-structures are well-understood locally and one of the main problem is to understand 1-flat G-structures which are not locally symmetric.
This problem was completely solved in \cite{MS} when the structure group $G \subset {\rm GL}(V)$ is an irreducible representation.

To facilitate our discussion, it is convenient to use the notion of cone structures.
When the structure group of a G-structure is a subgroup $G \subset {\rm GL}(V)$, a submanifold $Z \subset \BP V$ preserved by the $G$-action gives rise to a $Z$-isotrivial cone structure in the following sense.

\begin{definition}\label{d.isotrivial}
 Fix a vector space $V$ and a smooth projective variety $Z \subset \BP V$.
Let $Y$ be a complex manifold with $\dim Y = \dim V$.
\begin{itemize} \item[(i)] A  $Z$-{\em isotrivial cone structure} on $Y$ is a fiber subbundle $\sC \subset \BP TY$ such that the fiber $\sC_y \subset \BP T_y Y$ is projectively isomorphic to  $Z \subset \BP V$ for every $y \in Y$. \item[(ii)]
A $Z$-isotrivial cone structure determines a $G$-structure on $Y$ with $G = \Aut(\widehat{Z}) \subset {\rm GL}(V),$ which we call the G-structure {\em associated with } the $Z$-isotrivial cone structure. \item[(iii)]
A $Z$-isotrivial cone structure is {\em locally flat} (resp. {\em 1-flat}, resp. {\em locally symmetric}) if the associated G-structure  is locally  flat (resp. 1-flat, resp. locally symmetric).
\end{itemize} \end{definition}

In terms of Definition \ref{d.isotrivial}, we can reformulate one of the key results of \cite{MS} as follows.

\begin{theorem}\label{t.MS}
Let $Z \subset \BP V$ be the highest weight variety of an irreducible representation $G \subset {\rm GL}(V)$.
Then a $1$-flat $Z$-isotrivial cone structure is  locally symmetric,  unless $Z$ is one of the following varieties.
\begin{itemize}
\item[(1)]  A generalized minuscule variety of osculating rank 2 (in the sense of Definition \ref{d.minuscule}).
    \item[(2)] A subadjoint variety (in the sense of Definition \ref{d.subadjoint}).
      \end{itemize}
    \end{theorem}

This follows from Theorem 6.12 of \cite{MS}. In fact, the Lie algebra of $\Aut(\widehat{Z})$ is of the form $\C \oplus \fg$ for some semisimple Lie algebra $\fg$. If the structure is not locally symmetric, then $\dim K(\fg) >1$ in the notation of Lemma 2.3 of \cite{MS}. Thus by Theorem 6.12 \cite{MS}, the semisimple Lie algebra $\fg$ must be one of those in Table 6 or Table 7 of \cite{MS}. It is easy to check that the corresponding highest weight varieties are (1) and (2) of Theorem \ref{t.MS}.

In this paper, we propose to study Theorem \ref{t.MS} in a setting arising from algebraic geometry of minimal rational curves on uniruled projective manifolds.
Generalizing Definition \ref{d.isotrivial},  we define a (not necessarily isotrivial) cone structure on a manifold $Y$ as a closed subvariety $\sC \subset \BP TY$ dominant over $Y$.
A large class of interesting examples of cone structures arise from algebraic geometry, as VMRT-structures of minimal rational curves on uniruled projective manifolds (see \cite{HM99}, \cite{Hw01} and \cite{Hw12} for introductory surveys). Several examples of  such  VMRT-structures give isotrivial cone structures on suitable Zariski open subsets  of the projective manifolds, usually not on the whole projective manifolds.  So it is convenient to introduce the following terminology.

\begin{definition}\label{d.isovmrt} Let $Z \subset \BP V$ be as in Definition \ref{d.isotrivial}.
For a uniruled projective manifold $X$ and a family of minimal rational curves on $X$ (in the sense of Definition \ref{d.VMRT}), suppose that the VMRT $\sC_x \subset \BP T_x X$ is projectively isomorphic to $Z \subset \BP V$ for a general $x \in X$.  Then there exists a nonempty Zariski open subset $X^o \subset X$ where the family of VMRT's $\cup_{x \in X^o} \sC_x$ defines a $Z$-isotrivial cone structure on $X^o$. We call it a {\em $Z$-isotrivial VMRT structure}.\end{definition}

We want to understand when the situation of Theorem \ref{t.MS} can arise as a $Z$-isotrivial VMRT structure. The case (1) of Theorem \ref{t.MS} has been clarified in the VMRT setting by the following result of Ngaiming Mok, which is contained in the proof of Main Theorem in Section 2 of \cite{Mk}.

\begin{theorem}\label{t.Mok}
 Let  $Z \subset \BP V$  be a generalized minuscule variety of osculating rank 2 (equivalently, the highest weight variety $Z$ of the isotropy representation $V$ of an irreducible Hermitian symmetric space). Then any $Z$-isotrivial VMRT structure is locally flat.
   \end{theorem}

In this paper, we describe the case (2) of Theorem \ref{t.MS} in the VMRT setting in the following way when $\dim V \geq 5$.

 \begin{theorem}\label{t.main}
  Let $Z \subset \BP V$  be a subadjoint variety, excluding the case when $\dim V =4$ and $Z$ is the twisted cubic in $\BP^3$. Then any $Z$-isotrivial VMRT structure is locally flat.
   \end{theorem}

Note that both Theorem \ref{t.Mok} and Theorem \ref{t.main} do not assume the 1-flatness of the $Z$-isotrivial cone structure. So they answer, for the specified $Z \subset \BP V,$ the  following problem proposed in Question 1.8 of \cite{Hw10} .

\begin{question}\label{q.Hw10}
Let $Z \subset \BP V$ be a linearly nondegenerate smooth  projective variety. Is it true that any  $Z$-isotrivial VMRT structure is locally flat? \end{question}

An affirmative answer to Question \ref{q.Hw10} has been obtained when $Z$ is a smooth hypersurface (\cite{Hw10}, \cite{Hw13}) or a complete intersection of certain types (\cite{FH18}),
by the coframe method of \cite{Hw10}, which requires $\dim \aut(\widehat{Z}) =0.$
The problem becomes  much subtler when $\dim \aut(\widehat{Z}) >0,$ and there are many examples of homogeneous varieties $Z\subset \BP V$ for which the answer to Question \ref{q.Hw10} is negative.
Theorems \ref{t.Mok} and \ref{t.main} are the only known affirmative cases on Question \ref{q.Hw10} with $\dim \aut(\widehat{Z}) >0.$

As a corollary of Theorems \ref{t.MS}, \ref{t.Mok} and \ref{t.main}, we obtain

  \begin{corollary}\label{c.main} A 1-flat irreducible G-structure  that is given by varieties of minimal rational tangents  is locally symmetric, except possibly when the manifold has dimension $4$ and the highest weight variety is the twisted cubic. \end{corollary}

It should be mentioned that there are locally symmetric VMRT structures which are not locally flat:  the wonderful compactifications of simple algebraic groups have such VMRT structures by \cite{BF}. It is an interesting question to determine which locally symmetric structures in Theorem \ref{t.MS} do arise as VMRT structures.

It is natural to ask what one can say in the 4-dimensional case of Theorem \ref{t.main} and Corollary \ref{c.main}. We explain below why this case has been excluded from our discussion. It is still possible that the statements in Theorem \ref{t.main} or Corollary \ref{c.main}  hold including the 4-dimensional case.  Another natural question is whether Corollary \ref{c.main} can be generalized to a (not necessarily irreducible) 1-flat isotrivial cone structure given by varieties of minimal rational tangents. We leave them for future investigations.

Let us explain the ideas of the proof of Theorem \ref{t.main}.  We start recalling  the key components of the proof of Theorem \ref{t.Mok}.  Let $Z$ be as in Theorem \ref{t.Mok} and $\sC \subset \BP T X$ be the VMRT on a uniruled projective manifold with respect to a family $\sK$ of minimal rational curves such that the restriction $\sC|_{X^o}$ to a nonempty Zariski open subset $X^o \subset X$  is a $Z$-isotrivial cone structure.  Roughly speaking, the proof of Theorem \ref{t.Mok} in \cite{Mk} has the following two components,
 corresponding to the differential geometric part and the algebraic geometric part.

\begin{itemize} \item{{\sf Connection}.  A $Z$-isotrivial cone structure admits a holomorphic Cartan connection.}
\item{{\sf Extension}. The $Z$-isotrivial VMRT structure on $X^o$ (in the notation of Definition \ref{d.isovmrt}) can be extended to a $Z$-isotrivial cone structure on a neighborhood $U \subset X$ of a general member of $\sK$. }
 \end{itemize}

Once these two parts are established, we can show that the curvature of the Cartan connection restricted to $U$ vanishes by restricting it to the minimal rational curves in $U$, thereby prove the local flatness.  Can this argument work when $Z, X, X^o, \sK$ are  in the setting of Theorem \ref{t.main}? {\sf Connection} part is OK. In fact, such connections exist for $Z$-isotrivial cone structures for any highest weight variety $Z \subset \BP V$.
The problem is that {\sf Extension} part no longer works.  In fact, in an explicit example, like Example 1.7 of \cite{Hw10}, of the uniruled projective manifold $X$ with a $Z$-isotrivial VMRT structure for a highest weight variety $Z$, one can check that  unless $Z$ is the one in Theorem \ref{t.Mok},
    the $Z$-isotrivial cone structure always degenerates along a hypersurface in $X$ which intersects minimal rational curves nontrivially. So the $Z$-isotrivial cone structure cannot be extended to a neighborhood of a minimal rational curve in the setting of Theorem \ref{t.main}.

Our idea to overcome this difficulty is to {\em consider a geometric structure on an open subset $\sK^o \subset \sK$ in the space $\sK$ of minimal rational curves instead of the G-structure on  $X^o \subset X$}. The geometric structure on $\sK^o$, consisting of a distribution $D$ and an isotrivial cone structure $\sS \subset \BP D$, is related to the G-structure on $X^o$ via the universal family of minimal rational curves (see Section \ref{s.vmrt} for details).
$$ \begin{array}{ccccccccc}
& & & & {\rm Univ}_{\sK} & & & & \\
& & &\swarrow &    &\searrow& & &\\
(\sS \subset \BP D) \Rightarrow & \sK^o \subset &\sK& & & & X & \supset X^o & \Leftarrow \mbox{ G-strucure}
\end{array} $$
The idea of transferring a geometric structure on a complex manifold to a geometric structure of different type on the space of rational curves  goes back to Penrose transform in twistor theory. But it has not been used systematically in the study of minimal rational curves.
   In this respect, it is a novel feature of the current paper, compared with previous works like \cite{Mk} or \cite{HL}. The open subset $\sK^o$ contains a family of projective submanifolds biregular to the subadjoint variety.  Once we have a Cartan connection on $\sK^o$ for the geometric structure $\sS \subset \BP D$,  we can  argue as before by restricting the curvature tensor to the family of subadjoint varieties to show that the curvature vanishes. This means that the geometric structure  $\sS \subset \BP D$ is locally flat in a suitable sense, from which we deduce the local flatness of the VMRT structure on $X^o$.

But then the difficulty is transferred to  {\sf Connection} part. It is by no means clear that the new geometric structure $\sS \subset \BP D$ on $\sK^o$ admits a Cartan connection. Fortunately, we can apply  here the method developed  in  \cite{HL}, where for a  general class of geometric structures, the obstructions to the existence of Cartan connections are described as certain holomorphic tensors, under the assumption that  the underlying graded Lie algebras have finite-dimensional universal prolongations. We restrict these holomorphic tensors to the family of subadjoint varieties in $\sK^o$ to show that they vanish. The curvatures of the resulting connection vanish by the same argument.
The main reason that we have excluded the twisted cubic in $\BP^3$ in Theorem \ref{t.main} is that the corresponding graded Lie algebra has infinite-dimensional universal prolongation and the method of \cite{HL} cannot be applied.

As a side remark, we would like to mention that another approach to Theorem \ref{t.Mok}  avoiding {\sf Extension } part has been given in \cite{HN}, using the machinery of parabolic geometry. This approach cannot be used here because neither the G-structure in Theorem \ref{t.main}  nor the structure $\sS \subset \BP D$ on $\sK^o$ explained above are parabolic geometries.

The content of the paper is as follows.
In Section \ref{s.minuscule}, we review the basic properties of subadjoint varieties.
In Section \ref{s.fg}, we introduce the graded Lie algebra $\fg$, which is the underlying Lie algebra of  the geometric structure on $\sK^o,$ to be defined in Section \ref{s.vmrt}.   In Section \ref{s.HL}, we review the result in \cite{HL} of identifying obstructions to the existence of Cartan connections and curvatures for a general class of geometric structures.
In Section \ref{s.prolong}, we check the prolongation properties of the graded Lie algebra $\fg$, a condition necessary to apply the result in Section \ref{s.HL}.
 Section \ref{s.Spencer} identifies the components of the obstruction tensors that vanish when restricted to the subadjoint varieties in $ \sK^o$
and Section \ref{s.M} checks the vanishing of the remaining components by considering families of subadjoint varieties. In Section \ref{s.vmrt}, we explain how the setting of Theorem \ref{t.main} leads to the geometric structures studied in the previous sections. Note that rational curves and VMRT appear only in  Section \ref{s.vmrt}.

\section{Generalized minuscule varieties and subadjoint varieties}\label{s.minuscule}
In this section, we recall the definition and some basic properties of subadjoint varieties.
They are special cases of generalized minuscule varieties (i.e. Hermitian symmetric spaces equipped with equivariant projective embeddings). We start by introducing some notation related to generalized minuscule varieties.

\begin{notation}\label{n.hss}
Let $\fl = \fl_{-1} \oplus \fl_0 \oplus \fl_{1}$ be a graded semi-simple Lie algebra.
\begin{itemize} \item[(1)] Let $L$ be a  connected algebraic group with Lie algebra $\fl$ and let $P \subset L$ be the parabolic subgroup with Lie algebra $\fp= \fl_{-1} \oplus \fl_0$. The homogeneous space $Z= L/P$ is a Hermitian symmetric space. We can assume that $L$ acts effectively on $Z$ and $P$ is connected.  Denote by $L_{-1} \subset P$ the unipotent radical of $P$ with Lie algebra $\fl_{-1}$ and by $L_0 = P/L_{-1}$ the reductive group with Lie algebra $\fl_0$. Then the connected group $L_0$ acts effectively on the tangent space $T_z Z$ at the base point $z=[P] \in L/P.$
\item[(2)] Fix a set of simple roots $\{ \alpha_1, \ldots, \alpha_n\}$ and the fundamental weights $\{ \omega_1, \ldots, \omega_n\}$ of $\fl$ with respect to a fixed Cartan subalgebra. Let $\fl = \oplus_{i \in I} \fl^i$ be the decomposition into simple Lie algebras for a finite index set $I$.  For each $i \in I$, the parabolic subalgebra $ \fp^i= \fp \cap \fg^i$ is associated with some simple root $\alpha_i.$ Let $P^i \subset L^i$ be the corresponding parabolic subgroup of the simple Lie group $L^i$ and let $P_{\alpha_i} \subset L$ be the corresponding parabolic subgroup of $L$ such that $$L/P = \prod_{i \in I} L^i/P^i = \prod_{i\in I} L/P_{\alpha_i}.$$
\item[(3)] In terms of the decomposition in (2), there is a canonical identification ${\rm Pic}(L/P) = \oplus_{i\in I} \Z \omega_i$, where $\omega_i$ is the fundamental weight corresponding to the root $\alpha_i$.  A weight $\omega \in \oplus_{i\in I} \Z \omega_i$ determines a 1-dimensional $P$-module  which induces a line bundle  $\sO(\omega)$ on $L/P$. For $\omega = \sum_{i \in I} d_i \omega_i$, the line bundle $\sO(\omega)$  is very ample if and only if $d_i >0$ for all $i \in I$. In fact, the space of sections $H^0(L/P, \sO(\omega))$ is the irreducible $L$-module $V({\omega})$ with the highest weight $\omega$. This induces an $L$-equivariant embedding $L/P \subset \BP V({\omega})^{\vee}$.  Recall that $V({\omega})^{\vee} = V(- \bw_0(\omega))$ where $\bw_0$ is the longest element of the Weyl group of $L$. \end{itemize} \end{notation}

\begin{lemma}\label{l.xvv}
In Notation \ref{n.hss}, let $\sB \subset \BP \fl_1$ be the union of closed orbits under the $L_0$-action on $\BP\fl_1$ and let $\widehat{\sB} \subset \fl_1$ be its affine cone. Then
\begin{equation}\label{e.xvv} \{a \in \fl_{-1}, \ [[a, b], b] = 0 \mbox{ for all } b \in \widehat{\sB}\}  =  0. \end{equation}
\end{lemma}

\begin{proof}
Let us denote by $U $ the left-hand side of (\ref{e.xvv}), which
is a vector subspace of $\fl_{-1}$.  Since $L_0$ acts on $\fl$ preserving $\widehat{\sB} \subset \fl_1,$ we see that $U \subset \fl_{-1}$ is
an $L_0$-submodule.

Let $\fl= \oplus_{i \in I} \fl^i$ be the decomposition into simple Lie algebras, which induces the corresponding decompositions $$\fl_{-1} = \oplus_{i\in I} \fl_{-1}^i \mbox{ and } \fl_1 = \oplus_{i\in I} \fl_1^i.$$ Then $\sB \subset \BP \fl_1$ is the disjoint union of $\sB^i \subset \BP \fl_1^i, i \in I$, where $\sB^i$ is the unique closed orbit of the $L_0$-action on  $\BP \fl_1^i$.
An element $u \in U$ can be written as $u = \sum_{i \in I} a^i$ with $a^i \in \fl^i_{-1}$ such that $[[a^i, b^i], b^i] =0$ for all $b^i \in \widehat{\sB}^i$.
Thus it suffices to prove the proposition when $\fl$ is simple.

When $\fl$ is simple,  the subalgebra  $\fl_{-1}$ is an irreducible $\fl_0$-module and $\widehat{\sB}$ contains a positive root space $\fl_{\alpha}$ of $\fl$ such that $\fl_{-\alpha} \subset \fl_{-1}$ and the subalgebra   $$\fl_{-\alpha} \oplus \fl_{\alpha} \oplus [\fl_{-\alpha}, \fl_{\alpha}] \ \subset \ \fl$$
is isomorphic to $\fsl_2$, by the structure theory of  simple graded Lie algebras (e.g. Section 3.1 of \cite{Ya}). Then  $[[a,b], b] \neq 0$ for any nonzero $a \in \fl_{-\alpha}$
and nonzero $b \in \fl_\alpha$. It follows that  $\fl_{-\alpha} \cap U = 0$. Thus $U =0$ by the irreducibility of the $L_0$-module $\fl_{-1}.$
\end{proof}

The following definition is a slight modification of Definition 2.7 of \cite{LM03}.

\begin{definition}\label{d.minuscule}
In Notation \ref{n.hss} (3), fix a very ample line bundle $\sO(\omega)$   on $L/P$ and set $V = V(\omega)^{\vee}.$ The image $Z \subset \BP V$ of the embedding of $L/P$ by $\sO(\omega)$ is called a {\em generalized minuscule variety}.
   Fix a base point $z \in Z$ corresponding to a highest weight vector $v_0 \in V$ such that $P$ is the stabilizer of $z$ and $Z = L\cdot z$.   Denote by $|I|$ the cardinality of the index set $I$. \begin{itemize}
\item[(1)] The osculating spaces (see Section 2.2 of \cite{LM03}) of $Z \subset \BP V$ at $z$ are denoted by $$ 0 \subsetneq \widehat{z} \subsetneq \widehat{T}_z Z \subsetneq \widehat{T}_z^{(2)} Z \subsetneq \cdots \subsetneq \widehat{T}^{(r)}_z Z = V.$$
Denote by $$ V = V_0 \oplus V_1 \oplus \cdots \oplus V_{r}$$ the $L_0$-module decomposition corresponding to the osculating spaces. Let us call the positive integer $r$  the {\em osculating rank} of $Z \subset \BP V.$ \item[(2)] Denoting by $N^{(k)}_{Z,z}$ the $k$-th normal space  of $Z$ at $z$, we have natural identifications $$V_0 = \widehat{z} \mbox{ and } N^{(k)}_{Z,z}  = \Hom(V_0, V_k) \mbox{ for } 1 \leq k \leq r.$$
\item[(3)] The $k$-th projective fundamental form of $Z $ at $z$ is denoted by
$${\rm FF}^k_{Z,z}  \in \Hom ( \Sym^k T_z Z, N^{(k)}_{Z,z} )$$ for $2 \leq k \leq r$.
We also write ${\rm II}_{Z,z}$ for ${\rm FF}^2_{Z,z}$ and ${\rm III}_{Z,z}$ for ${\rm FF}^3_{Z,z}$.
\item[(4)] When $L$ is simple and $|I|=1$, the image $Z \subset \BP V(\omega_i)^{\vee}$  of the embedding of $L/P= L/P_{\alpha_i}$ by the ample generator $\sO(\omega_i)$  of ${\rm Pic}(L/P)$ is   called a {\em minuscule variety}. \end{itemize}   \end{definition}

We are particularly interested in the following class of generalized minuscule varieties.

\begin{definition}\label{d.subadjoint}
Each simple Lie algebra $\fs$ admits a unique gradation (called the contact gradation in Section 4.2 of \cite{Ya}) of the form $$\fs = \fs_{-2} \oplus \fs_{-1} \oplus \fs_0 \oplus \fs_1 \oplus \fs_2 \mbox{ with } \dim \fs_2 = \dim \fs_{-2} =1.$$
The Lie bracket $\fs_{1} \wedge \fs_{1} \to \fs_{2}$ is a  symplectic form on $\fs_{1}$. Let us exclude $\fs$  of type $A$ or $C$ from our discussion. Then the natural representation of $\fs_0$ on $\fs_{1}$ is irreducible and its highest weight variety $Z \subset \BP \fs_{1}$  is a generalized minuscule variety of osculating rank 3. The projective variety $Z \subset \BP \fs_{1}$ is called a \emph{subadjoint variety}.  \end{definition}

The following proposition is well-known (e.g. Section 1.4.6 of \cite{Hw01} or Table 1 in \cite{Bu}).

\begin{proposition}\label{p.list}
The explicit description of subadjoint varieties is as follows. Here, by abuse of notation, we denote a simple Lie group by the type of its Dynkin diagram.

\begin{itemize}
    \item[(0)] When $\fs$ is of type $G_2$, the subadjoint variety $Z \subset \BP \fs_{1}$ is isomorphic to $A_1/P_{\alpha_1} \subset \BP V(3 \omega_1)^{\vee}$, which is the twisted cubic curve in $\BP^3$.
        \item[(i)] The following four subadjoint varieties are minuscule varieties.
         \begin{itemize} \item[($F_4$)] When $\fs$ is of type $F_4$, the subadjoint variety $Z \subset \BP \fs_{1}$ is isomorphic to the minuscule variety  $C_3/P_{\alpha_3} \subset \BP V(\omega_3)^{\vee}$, which is  the Pl\"ucker embedding of the Lagrangian Grassmannian of 3-dimensional isotropic subspaces in  a symplectic vector space of dimension 6.
\item[($E_6$)] When $\fs$ is of type $E_6$, the subadjoint variety $Z \subset \BP \fs_{1}$ is isomorphic to the minuscule variety $A_5/P_{\alpha_3} \subset \BP V(\omega_3)^{\vee}$, which is  the Pl\"ucker embedding of the Grassmannian  of 3-dimensional subspaces in a  vector space of dimension 6.
    \item[($E_7$)] When $\fs$ is of type $E_7$, the subadjoint variety $Z\subset \BP \fs_{1}$ is isomorphic to the minuscule variety $D_6/P_{\alpha_6} \subset \BP V(\omega_6)^{\vee},$ which is the  embedding of the spinor variety of
 isotropic subspaces of dimension 6 in an orthogonal vector space of dimension 12.
\item[($E_8$)] When $\fs$ is of type $E_8$, the subadjoint variety $Z\subset \BP \fs_{1}$ is isomorphic to the 27-dimensional minuscule variety $$E_7/P_{\alpha_7} \subset \BP V(\omega_7)^{\vee} \cong \BP^{55}.$$
\end{itemize} \item[(ii)] When $\fs$ is $\fs \fo_{m+4}$ with $m \geq 3$, the subadjoint variety $Z \subset \BP \fs_{1}$ is isomorphic to  the Segre product $\BP^1 \times \Q^{m-2} \subset \BP^{2m-1}$.
\begin{itemize} \item[(ii-1)] When $m=3$, the surface $Z \subset \BP^5$ is isomorphic to the Segre product of a line and a conic, biregular to $\BP^1 \times \BP^1$.
\item[(ii-2)] When $m=4$, the threefold $Z \subset \BP^{2m-1}$ is the Segre product of three lines, biregular to $\BP^1 \times \BP^1 \times \BP^1.$
    \item[(ii-3)] When $m \geq 5$, the Segre product $Z \cong \BP^1 \times \Q^{m-2}$ is isomorphic to  the generalized minuscule variety  $$A_1/P_{\alpha_1} \times B_{\ell}/P_{\alpha_1} \subset \BP V(\omega^A_1 + \omega^B_1)^{\vee} \mbox{ for } m = 2 \ell +1, \mbox{ or } $$  $$A_1/P_{\alpha_1} \times D_{\ell}/P_{\alpha_1} \subset \BP V(\omega^A_1 + \omega^D_1)^{\vee} \mbox{ for } m= 2 \ell,$$ where $\omega^A_1, \omega^B_1, \omega^D_1$ denote the first fundamental weights of the Lie algebras of type $A, B, D$, respectively.
\end{itemize}
\end{itemize} \end{proposition}

The following properties of the fundamental forms of subadjoint varieties are well-known.

\begin{proposition}\label{p.FFsubadjoint}
Let $Z \subset \BP \fs_{1}$ be a subadjoint variety,  excluding the case (0) in Proposition \ref{p.list}.
Let $\sigma: \wedge^2 \fs_1 \to \fs_2$ be the symplectic form in Definition \ref{d.subadjoint}.
For each point $z \in Z$ and the corresponding 1-dimensional subspace $\widehat{z} \subset \fs_1$,  the following holds.
\begin{itemize} \item[(1)] The affine cone $\widehat{Z} \subset \fs_1$ is Lagrangian with respect to $\sigma$, i.e. it has dimension $\frac{1}{2} \dim \fs_1$ and $\sigma(\widehat{T}_{z}Z, \widehat{T}_{z}Z) = 0$ for the affine tangent space $\widehat{T}_{z}Z \subset \fs_1$ of $Z$ at $z$.  We can also say that $Z$ is a Legendrian submanifold of $\BP \fs_1$ equipped with the contact structure induced by $\sigma$.
\item[(2)] The second osculating space $T^{(2)}_z Z \subset \fs_{1}$ is the hyperplane annihilated by $\sigma(\widehat{z}, \cdot)$.
     \item[(3)] The third fundamental form ${\rm III}_{Z,z},$ which is a  cubic form on $T_z Z$ by (2), is
     nondegenerate in the sense that
        $$\{a \in T_z Z \mid  \III_{Z,z} (a, b, c) = 0 \mbox{ for all } b, c \in T_{z}Z \} = 0.$$
    \item[(4)] The symplectic form $\sigma$ determines a perfect pairing $\beta: T_{z}Z \otimes N^{(2)}_{Z,z} \to N^{(3)}_{Z,z} \cong \C$ such that for any $a_1, a_2, a_3 \in T_z Z$, we have
    $$\beta(a_1, \II_{Z,z}(a_2, a_3) ) =  \III_{Z,z}(a_1, a_2, a_3).$$
    Consequently, if there exists a bilinear form $\beta': T_{z}Z \otimes N^{(2)}_{Z,z} \to \C$ such that for some nonzero $v \in N^{(3)}_{Z,z}$,
     $$\beta'(a_1, \II_{Z,z}(a_2, a_3) ) v =  \III_{Z,z}(a_1, a_2, a_3) \mbox{ for all } a_1, a_2, a_3 \in T_{z}Z,$$ then $\beta'$ is a perfect pairing.
  \item[(5)] When $Z = L/P$ as a generalized minuscule variety, let $\sB_z \subset \BP T_z Z$ be the union of all closed orbits of the action of $P $ on $\BP T_z Z$. Define  the base loci ${\rm Bs}(\II_{Z,z}), {\rm Bs}(\III_{Z,z}) \subset \BP T_z Z$  of the second and third fundamental forms as    \begin{eqnarray*} {\rm Bs}(\II_{Z,z}) & = & \BP \{ u \in T_z Z \mid \II_{Z,z}(u, u) = 0 \} \\ {\rm Bs}(\III_{Z,z}) &=& \BP \{ u \in T_z Z  \mid \III_{Z,z}(u,u, u) =0\}. \end{eqnarray*}  Then $ {\rm Bs}(\II_{Z,z}) \subsetneq\sB_z = {\rm Bs}(\III_{Z,z})$ in the case (ii-1) of Proposition \ref{p.list} and ${\rm Bs}(\II_{Z,z})  = \sB_z$ is the singular locus of the cubic hypersurface ${\rm Bs}(\III_{Z,z})$ in $\BP T_z Z$ in the other cases of Proposition \ref{p.list}.  
      \end{itemize}  \end{proposition}

      \begin{proof}
      (1) can be found in Section 1.4.6 of \cite{Hw01} or Table 1 in \cite{Bu}.
       (2), (3) and (4) can be checked easily in each case of Proposition \ref{p.list}. They are also properties of projective Legendrian submanifolds with nondegenerate second fundamental forms and can be deduced from pp. 348-349 of \cite{LM07} (see also Proposition 2 in \cite{Hw21}).   (5) can be checked by elementary arguments for the case (ii) of Proposition \ref{p.list}.
        For the four cases in (i) of Proposition \ref{p.list},  the four projective varieties $\sB_z \subset \BP T_z Z$ are exactly the Severi varieties. Recall that the secant variety of a Severi variety is a cubic hypersurface whose singular locus is exactly the Severi variety (e.g. Theorem 2.4 in Chapter IV of \cite{Za}).    The base locus of the third fundamental form is the unique cubic hypersurface invariant under the isotropy action and must be isomorphic to  the secant  variety of the Severi variety. 
       \end{proof}

\section{Graded Lie algebra $\fg$ associated to a subadjoint variety}\label{s.fg}

The following graded Lie algebra plays a key role in our discussion.

  \begin{definition}\label{d.fh}
        Let us denote by $\bV$ the subspace $\fs_{1}$  in Definition \ref{d.subadjoint} and by $Z  \subset \BP \bV$ the subadjoint variety. Throughout, we exclude the case (0) of Proposition \ref{p.list} in our discussion. Fixing a lowest weight vector $v_0 \in \bV$, we can write $Z = L/P$ with  the associated graded semisimple Lie algebra $\fl= \fl_{-1} \oplus \fl_0 \oplus \fl_{1}.$
Let $\bV = \bV_0 \oplus \bV_{1} \oplus \bV_{2} \oplus \bV_{3}$ be the decomposition into irreducible $\fl_0$-modules.
        Let $\C {\rm Id}_{\bV} \oplus \fl \subset \fgl(\bV)$ be the Lie algebra of the group $\Aut(\widehat{Z}) \subset {\rm GL}(\V)$ of linear automorphisms preserving $\widehat{Z} \subset \bV$.
The Lie algebra $\fg := (\C {\rm Id}_{\bV} \oplus \fl) \ltimes \bV$ has the following graded Lie algebra structure:
        \begin{eqnarray*} \fg_{-1} & = & \fl_{-1} \\ \fg_0 &=& \bV_0 \oplus \C {\rm Id}_{\bV} \oplus \fl_0 \\
        \fg_{1} & =& \bV_{1} \oplus \fl_{1} \ = \ (\fl_{1}  \otimes \bV_0) \oplus \fl_{1}\\ \fg_{2} & = & \bV_{2} \\ \fg_{3} &=& \bV_{3}. \end{eqnarray*}  Let us write $ a_1 \cdot v_0 \in \bV_1$, $a_1a_2 \cdot v_0 \in \bV_2$ and $a_1a_2a_3 \cdot v_0 \in \bV_3$ for $[a_1, v_0]$, $[a_1, [a_2, v_0]]$ and $[a_1, [a_2, [a_3, v_0]]]$ respectively, where $a_1, a_2, a_3 \in \fl_1.$ Denote by $\bc: \fl_0 \to \C$ the linear functional satisfying $$[ b, v_o] = \bc(b) v_o \mbox{ for all } b \in \fl_0.$$
                \end{definition}

                \begin{remark}\label{remark1}
     The motivation for considering the graded Lie algebra $\fg$ in Definition \ref{d.fh} is the following.
  Consider the flat $Z$-isotrivial cone structure $\sC^{\rm flat} \subset \BP \bV$ discussed in Section \ref{s.I} associated to our subadjoint variety $Z \subset \BP \bV$. The group $\bG= (\C^{\times} \times L) \ltimes \bV$ acts naturally on the manifolds $\bV$ and $\sC^{\rm flat}$. Let $\bK$ be the space of affine lines on $\bV$ in the direction of $\widehat{Z} \subset \bV$.
  The tangent spaces of such affine lines lie in $\sC^{\rm flat}$ and determine a $\bG$-equivariant fibration $\sC^{\rm flat} \to \bK$. Our interest lies in the $\bG$-homogeneous space $\bK$. The Lie algebra of $\bG$ is $\fg$ in Definition \ref{d.fh} and the Lie subalgebra of the isotropy group at a base point of $\bK$ is $\fg_{-1} \oplus \fg_0$.
          \end{remark}

                \begin{lemma}\label{l.II}
        In Definition \ref{d.fh}, we have  homomorphisms ${\rm II}: \Sym^2 \fl_{1} \to \bV_0^{\vee} \otimes \bV_{2}$ and ${\rm III}: \Sym^3 \fl_1 \to \bV_0^{\vee} \otimes \bV_3$ such that the Lie bracket in $\fg$ satisfies
       $$ [a, b \cdot v_0 ] = {\rm II}(a, b)\cdot v_0 $$ $$  [a, [ b, c \cdot v_0]]
      =  {\rm III}(a,b,c) \cdot v_0 $$
for all $a, b, c \in \fl_{1},$ where the dot on the right hand side stands for the contraction of $v_0$ with elements of $\bV_0^{\vee}$.   Moreover, \begin{equation}\label{e.II} [a+ t a \cdot v_0, b+ t b \cdot v_0] =0 \end{equation} for any $a, b \in \fl_1$  and $t \in \C$.  \end{lemma}

        \begin{proof} The second fundamental form (resp. the third fundamental form) of $Z$ at $[v_0]$ induces a homomorphism $${\rm II}: \Sym^2 \fl_{-1} \to \bV_0^{\vee} \otimes \bV_{2} \mbox{ (resp. } {\rm III}: \Sym^3 \fl_{-1} \to \bV_0^{\vee}\otimes \bV_3 {\rm )}.$$
        The relation of the Lie bracket in $\fg$ with ${\rm II}$ and ${\rm III}$ follows from the relation of the representation of $\fl$ on $\bV$ and the fundamental forms of $Z = L/P \subset \BP \bV$ given in Proposition 2.3 of \cite{LM03}.  The last statement follows from
\begin{eqnarray*}
[a+ t a \cdot v_0, b+ tb \cdot v_0] &=& t [a, b \cdot v_0] - t[b, a \cdot v_0]\\
&=& t{\rm II}(a,b)\cdot v_o - t {\rm II}(b,a)\cdot v_0 \  = \ 0.
\end{eqnarray*}
\end{proof}

        \begin{proposition}\label{p.auttimes}
       In Definition \ref{d.fh},  let $\sB \subset \BP \fl_{1}$ be the union of the highest weight varieties of the representation of $\fl_0$ on $\fl_{1}$.
        Denote by $\fg_+$ the nilpotent graded Lie algebra $\fg_{1} + \fg_{2} + \fg_{3}$ and by ${\rm gr}\Aut(\fg_+)$ the graded Lie algebra automorphism group of $\fg_+$.
        \begin{itemize}
        \item[(1)] The Lie algebra of $\Aut(\widehat{\sB}) \subset {\rm GL}(\fl_{1})$ is isomorphic to $\fl_0$ via the adjoint representation of $\fl$.
        \item[(2)] The restriction homomorphism ${\rm gr}\Aut(\fg_+) \to {\rm GL}(\fg_{1})$ is injective.
        \item[(3)] Setting $\bW := \bV_0 \oplus \C$, we have  a tensor decomposition $\fg_{1} = \bW \otimes \fl_{1}$  and the Segre embedding $\BP \bW \times \BP \fl_1 \subset \BP \fg_{1}.$ Define the subgroup ${\bG}_0 \subset {\rm gr}\Aut(\fg_+)$ as the connected component of the group $$ \widetilde{\bG}_0:= \{ g \in {\rm gr} \Aut(\fg_+) \mid \   g (\BP \bW \times \BP \fl_1) = \BP \bW \times \BP \fl_1\}.$$
            Then $\fg_0$ is the Lie algebra of ${\bG}_0$.
         \end{itemize}
         \end{proposition}

         \begin{proof}
         (1) is easy to check from the list  of subadjoint varieties in Proposition \ref{p.list}. (2) is immediate from the fact that $\fg_{1}$ generates $\fg_+$ as a Lie algebra.

         To see (3), consider the natural homomorphism $\Phi: {\rm GL}(\bW) \times {\rm GL}(\fl_{1}) \to {\rm GL}(\fg_{1})$, which has a 1-dimensional kernel. The subgroup of ${\rm GL}(\fg_{1})$ preserving $\BP \bW \times \BP \fl_1$ is $\Phi({\rm GL}(\bW) \times {\rm GL}(\fl_1))$.     Define $${\rm GL}(\bW)_{[\bV_0]} := \{ g \in {\rm GL}(\bW) \mid  g(\bV_0) = \bV_0\}.$$

         We claim that $$\bV_1 = \{ u \in \fg_1 \mid [u, \bV_2] =0\}.$$
         Otherwise, there exists $0 \neq a \in \fl_1$ satisfying $[a, \bV_2] =0.$
         Then for any $ b, c \in \fl_1$, we have $$[a, [b, c \cdot v_0]] = \III(a, b, c) \cdot v_0 =0,$$ a contradiction to the nondegeneracy of $\III$ in Proposition \ref{p.FFsubadjoint} (3).

By the claim, the subspace $\bV_1$ is preserved by the action of the group ${\rm gr}\Aut(\fg_+)$. It follows that the image of $G_0$ in $ {\rm GL}(\fg_1)$ is contained in $\Phi({\rm GL}(\bW)_{[\bV_0]} \times {\rm GL}(\fl_{1})).$

By Lemma \ref{l.II}, the image of ${\rm gr}\Aut(\fg_+)$ in ${\rm GL}(\fg_1)$ preserves $\BP\bW\times{\rm Bs}(\II)$ and $\BP\bW\times{\rm Bs}(\III)$. By (1) and Proposition \ref{p.FFsubadjoint}, this implies that
the image of $\bG_0$ is contained in $\Phi({\rm GL}(\bW)_{[\bV_0]} \times L_0),$ which has dimension $2 + \dim \fl_0 = \dim \fg_0$.
         It is easy to see that $\fg_0$ is a subalgebra of the Lie algebra of ${\bG}_0$ via the adjoint representation of $\fg_0$ on $\fg_+.$
         As $\dim {\bG}_0 \leq \dim \fg_0$, we see that $\fg_0$ is the Lie algebra of ${\bf G}_0$. \end{proof}

\section{Review of Section 2 of \cite{HL}}\label{s.HL}

As explained in Remark \ref{remark1}, we have a certain homogeneous space $\bK$ associated with the Lie algebra $\fg$ in Definition \ref{d.fh}. This homogenous space has a natural geometric structure
(see Example \ref{ex.model}).  The proof of our main result Theorem \ref{t.main} is eventually reduced to check that a manifold having a geometric structure analogous to that of $\bK$ is actually locally isomorphic to $\bK$, when the manifold contains a family of submanifolds of a certain type. To check such a statement, we employ the method developed in Section 2 in \cite{HL}, which we recall  in this section.
Note that we use a  sign convention for the indices of graded Lie algebras, which is  different from the one used in \cite{HL} or \cite{Ya},  because this makes the notation much simpler in the setting of this paper.
It is straightforward to translate the discussion in Section 2 of \cite{HL} into the following discussion, just by switching the signs of the indices.

\begin{notation}\label{n.grade}
When we work with a graded vector space $V = \oplus_{i=a}^{i=b}  V_i$ with $a, b \in \Z$, we sometimes write it as $V = \oplus_{i \in \Z} V_i$ with the understanding that $V_i=0$ if $i <a$ or $i >b$.
For two graded vector spaces $V = \oplus_{i\in \Z} V_i$ and $W = \oplus_{j \in \Z} W_j$, and an integer $k$, denote by $\Hom(V, W)_k$ the vector space of homomorphisms of degree $k$ from $V$ to $W$.
\end{notation}

\begin{definition}\label{d.prolong}
Let $d$ be a positive integer and let $\fn_{+} = \oplus_{i=1}^{d} \fn_i$ be a graded nilpotent Lie algebra.
Denote by ${\rm gr}\Aut(\fn_+)$ the group of graded Lie algebra automorphisms of $\fn_+$ and  by ${\rm gr}\aut(\fn_+)$ its Lie algebra.  For a subalgebra $\fn_0 \subset {\rm gr}\aut(\fn_+)$, the {\em universal prolongation} of $\fn_0 \oplus \fn_+$ is the graded Lie algebra  $$\fn = \bigoplus_{k= 1}^{\infty}  \fn_{-k} \oplus \fn_0 \oplus \fn_+ $$  such that for each $k >0$,
the Lie bracket determines an injection $\fn_{-k} \subset \Hom(\fn_+, \fn)_{-k}$ whose image is the space of all linear homomorphisms of degree $-k$  satisfying \begin{equation}\label{e.prolong} \varphi([x, y]) = [\varphi(x), y] + [x, \varphi(y)] \mbox{ for all } x, y \in \fn_+.\end{equation} The universal prolongation is unique up to isomorphism. The subspace $\fn_{-k}$ is called the $k$-th {\em prolongation} of $\fn_0 \oplus \fn_+.$  When $\fn_+$ is generated by $\fn_1$, the vanishing $\fn_{-k}=0$ implies that $\fn_{-i} =0$ for $i\geq k$. \end{definition}

\begin{definition}\label{d.filtration}
Let $M$ be a complex manifold. A {\em filtration } $F^{\bullet}$ on $M$ is a collection of subbundles
$$0= F^0 \subset F^1 \subset F^2 \subset \cdots \subset F^{d-1} \subset F^d = TM$$ of the tangent bundle $TM$ for some $d \in \N$ such that
viewing the locally free sheaf $\sO(F^i)$ as a sheaf of vector fields on $M$, we have \begin{equation}\label{e.filt} [\sO(F^{i}), \sO( F^{j})] \subset \sO(F^{k})\end{equation}  for the Lie bracket of vector fields and for any nonnegative integers
     $i, j, k,$
 satisfying $i+j \leq k$.
  For a point $x \in M$, denote by $F^{k}_x $ the fiber of $F^k $ at the point $x$. By (\ref{e.filt}), the graded vector space $${\rm Symb}_x(F^{\bullet}) \ := \ \bigoplus_{i \in \N} F_x^{i}/F_x^{i-1} $$
 for each point $x \in M$ has the structure of a nilpotent graded Lie algebra, called the {\em symbol algebra} of the filtration at $x$.
 \end{definition}

 \begin{definition}\label{d.G_0-structure}
For a graded nilpotent Lie algebra $\fg_{+}= \fg_{1} \oplus \cdots \oplus \fg_{d},$
a {\em filtration of type} $\fg_{+}$ on a complex manifold $M$ is a filtration $F^{\bullet}$  on $M$ as in Definition \ref{d.filtration} such that
 for any $x\in M$, the symbol algebra ${\rm Symb}_x(F^{\bullet})$  is isomorphic to $\fg_{+}$ as graded Lie algebras.
The {\em graded frame bundle} of the manifold $M$ with a filtration $F^{\bullet}$ of type $\fg_{+}$  is the  ${\rm gr}\Aut(\fg_{+})$-principal bundle ${\rm grFr}(F^{\bullet})$ on $M$ whose fiber at $x$ is the set of graded Lie algebra isomorphisms from
$\fg_{+}$ to ${\rm Symb}_x(F^{\bullet})$. Let $G_0 \subset {\rm gr}\Aut(\fg_{+})$ be a connected closed subgroup.  A $G_0$-{\em structure subordinate to the filtration} $F^{\bullet}$ on $M$ means a $G_0$-principal subbundle  $\sA \subset {\rm grFr}(F^{\bullet})$.\end{definition}

\begin{example}\label{ex.model}
In Definition \ref{d.G_0-structure}, assume that the universal prolongation $\fg:= \oplus_{i \leq d} \fg_i$ is finite-dimensional.
Let $G$ be a connected Lie group with Lie algebra $\fg$ and assume that the connected subgroup $G^0 \subset G$  with Lie algebra $\oplus_{i \leq 0} \fg_i$ is closed.
Then the homogenous space $G/G^0$ has a natural filtration $F^{\bullet}_{G/G^0}$ of type $\fg_+$ and a natural $G_0$-structure $\sA_{G/G^0} \subset {\rm grFr}(F^{\bullet}_{G/G^0}) $ subordinate to the filtration. \end{example}

\begin{definition}\label{d.cohomology}
Let $\fg= \oplus_{i \leq d} \fg_i$  be a graded Lie algebra which is the universal prolongation of $\fg_0 \oplus \fg_+$ in the sense of Definition \ref{d.prolong}.
Define, for each negative integer $k$,
\begin{eqnarray*}
C^{k, 1}(\fg) & := & \Hom(\fg_+, \fg)_k , \\
C^{k, 2}(\fg) & :=&
\Hom(\wedge^2 \fg_{+}, \fg)_k.
\end{eqnarray*}
For an element $f \in C^{k, 1}(\fg)$, define $\partial f \in C^{k, 2}(\fg)$ by
$$\partial f(u, v)  = [f(u), v] + [u, f(v)] -f([u, v]) \mbox{ for all } u, v \in \fg_{+}.$$ This determines a $G_0$-module homomorphism $$\partial: C^{k, 1}(\fg) \to C^{k, 2}(\fg).$$  Its cokernel $Q^k( \fg):= C^{k, 2}(\fg)/\partial(C^{k,1}(\fg))$ has a natural $G_0$-module structure. \end{definition}

The following is the main result in Section 2 of \cite{HL}: it is Theorem 2.17 combined with Proposition 2.5 in \cite{HL}.

\begin{theorem}\label{t.HL}
In Definition \ref{d.G_0-structure}, assume that the action of $G_0$ on $\fg_+$ has no nonzero fixed vector and the universal prolongation $\fg:= \oplus_{i \leq d} \fg_i$ is finite-dimensional. Let $G$ be a connected Lie group of adjoint type with Lie algebra $\fg$ and $G^0 \subset G$ be the connected closed subgroup with Lie algebra $\oplus_{i \leq 0} \fg_i$. Let $M$ be a simply connected complex manifold  with a filtration $F^{\bullet}$ of type $\fg_{+}$ and let $\sA \subset {\rm grFr}(F^{\bullet})$ be a $G_0$-structure on $M$ subordinate to the filtration. For each negative integer $k$, let $\sE_M(Q^k(\fg))$ be the vector bundle associated to the principal $G_0$-bundle $\sA$ by the natural representation of $G_0$ on
$Q^k(\fg)$ from Definition \ref{d.cohomology}.  Assume that $H^0(M, \sE_M(Q^k(\fg))) =0$ for all $ k \leq -1.$ Then there exists an open immersion $\phi: M \to G/G^0$ that sends $F^{\bullet}$ (resp.  $\sA$) isomorphically to $ F^{\bullet}_{G/G^0}|_{\phi(M)}$ (resp. $\sA_{G/G^0}|_{\phi(M)}$) in Example \ref{ex.model}.
\end{theorem}

The statement of Theorem \ref{t.HL} is slightly different from Theorem 2.17 of \cite{HL}. The latter theorem gives the existence of a Cartan connection on a prolonged principal bundle under the weaker assumption of the vanishing $$H^0(M, \sE_M(Q^k(\fg)) ) =0 \mbox{ for all } -m-d \leq k \leq - 1,$$ for some positive integer $m$.  By assuming the vanishing for all $k \leq -1$,  we can obtain the local flatness of the Cartan connection by the same argument as in the proof of Theorem 2.17 of \cite{HL}, because the curvature of the resulting Cartan connection would give elements of $H^0(M, \sE_M(Q^k(\fg)) )$ for some values of $k \leq -1$. This implies the existence of the open immersion $f$ by Proposition 2.5 of \cite{HL}.

        \section{Prolongation property of the graded Lie algebra $\fg$}\label{s.prolong}

The main result of this section is the following.

\begin{theorem}\label{t.prolong}
Let $\fg$ be the graded Lie algebra  in Definition \ref{d.fh}. Denote by $\fp_{-k}$ the $k$-th prolongation of $\fg_0 \oplus \fg_{+}$ for $k \geq 1$ in the sense of Definition \ref{d.prolong}. Then $\fp_{-1} \cong \fg_{-1}$ and $\fp_{-2} = 0$. In other words, the graded Lie algebra $\fg = \fg_{-1} \oplus \fg_0 \oplus \fg_+$ is the universal prolongation of $\fg_0 \oplus \fg_{+}$.
\end{theorem}

Before discussing the proof of Theorem \ref{t.prolong}, we recall some  well-known results on
universal prolongations.  The first is the following easy lemma on the direct sum of universal prolongations.

\begin{lemma}\label{l.sumprolong}
Let $\fq'_{1}$ and $\fq''_{1}$ be two vector spaces which we regard as abelian Lie algebras. Choose Lie subalgebras $\fq'_0 \subset \fgl(\fq'_{1})$ and $\fq''_0 \subset \fgl(\fq''_{1})$. Let $\fq'=  (\bigoplus_{k=1}^{\infty} \fq'_{-k}) \oplus \fq'_0 \oplus  \fq'_1$ (resp. $\fq'' = (\bigoplus_{k=1}^{\infty}\fq''_{-k}) \oplus  \fq''_0 \oplus  \fq''_1$) be the universal prolongation of $\fq'_{0} \oplus  \fq'_1$ (resp. $\fq''_{0} \oplus  \fq''_1$). Then for the  direct products of Lie algebras  $$\fq_{1} = \fq'_{1} \oplus \fq''_{1} \mbox{ and  } \fq_0 = \fq'_0 \oplus \fq''_0 \subset {\rm gr} \aut(\fq_1),$$ the universal prolongation $\fq = (\bigoplus_{k=1}^{\infty} \fq_{-k}) \oplus \fq_0 \oplus \fq_1 $ of $\fq_{0} \oplus \fq_{1}$
 is isomorphic to $\fq'  \oplus \fq''$. In other words, we have a natural identification $\fq_{-k} = \fq'_{-k} \oplus \fq''_{-k}$ for each $k \geq 1$.
\end{lemma}

\begin{proof}
We prove $\fq_{-i} = \fq'_{-i} \oplus \fq''_{-i}$ for all $i \geq 0$ by induction. There is nothing to prove for $i= 0$.
Assume that this has been checked for any $i\leq k-1$ for some $k \geq 1$, let us check it for $i= k$.
The induction hypothesis says that $\fq_{1-k} = \fq'_{1-k} \oplus \fq''_{1-k}$ which implies $$ \Hom(\fq_1, \fq)_{-k} = \Hom(\fq_1, \fq' \oplus \fq'')_{-k}. $$
Thus any  $$\varphi \in \fq_{-k} \subset \Hom(\fq_1, \fq)_{-k} = \Hom ( \fq'_{1} \oplus \fq''_{1}, \fq' \oplus \fq'')_{-k}$$
has a unique decomposition $$\varphi = \varphi' + \widehat{\varphi}' + \varphi'' + \widehat{\varphi}'',$$
such that $$\varphi' \in \Hom(\fq'_{1}, \fq')_{-k},  \ \varphi'' \in \Hom(\fq''_{1}, \fq'')_{-k} $$
$$\widehat{\varphi}' \in \Hom(\fq'_{1}, \fq'')_{-k}, \ \widehat{\varphi}'' \in \Hom(\fq''_{1}, \fq')_{-k}.$$
For any $a \in \fq'_{1}$ and $b \in \fq''_{1}$,
$$[\varphi'(a), b]= 0 = [a, \varphi''(b)].$$ Thus (\ref{e.prolong}) gives
\begin{eqnarray*} 0 = \varphi([a,b]) &=& [\varphi(a), b] + [a, \varphi(b)]  \\
& =& [\widehat{\varphi}'(a), b] + [a, \widehat{\varphi}''(b)], \end{eqnarray*}
which implies $ [\widehat{\varphi}'(a), b] =0 =  [a, \widehat{\varphi}''(b)].$
Since this holds for all $a \in \fq'_{-1}$ and $b \in \fq''_{-1}$, we conclude  $\widehat{\varphi}' = \widehat{\varphi}'' =0$ and $\varphi = \varphi'+ \varphi''$, which says $\fq_{-k} = \fq'_{-k} \oplus \fq''_{-k}$. \end{proof}

Next, we have the following special case of Theorem 5.2 of   \cite{Ya}.

\begin{theorem}\label{t.Yamaguchi}
Let   $\fl_{-1} \oplus \fl_0 \oplus \fl_{1}$ be the grading on a simple Lie algebra $\fl$ associated to a minuscule variety, different from projective spaces.  Then $\fl_{-1} \oplus \fl_0 \oplus \fl_{1}$ is the universal prolongation of $\fl_{0} \oplus \fl_1$. \end{theorem}

It is necessary to exclude projective spaces in Theorem \ref{t.Yamaguchi}. For example,
when the minuscule variety is $\BP^1$, we have the following result, which can be checked by a straightforward computation.

\begin{lemma}\label{l.sl2}
Let $\ff=  \bigoplus_{k\geq-1} \ff_{-k}$ be the graded Lie algebra of formal vector fields in one variable: each $\ff_{-k}$ is a 1-dimensional vector space with a basis $t^{k+1}\frac{\p}{\p t}$  for a formal variable $t$ and the Lie brackets satisfy
$$ [t\frac{\p}{\p t}, t^{k+1} \frac{\p}{\p t}] = k t^{k+1}\frac{\p}{\p t}, \ [t^2 \frac{\p}{\p t}, \frac{\p}{\p t}] = - 2 t \frac{\p}{\p t}, \ [\frac{\p}{\p t}, t^{k+1} \frac{\p}{\p t}] = (k+1) t^k \frac{\p}{\p t} $$ for all $ k \geq -1.$   Then the following holds.
\begin{itemize} \item[(1)] The Lie algebra
$\ff$ is the universal prolongation of $\ff_{0} \oplus \ff_1$.
\item[(2)] If $b \in \ff_{-k}, k \geq 1,$ satisfies $[[b, a], a] =0$ for a nonzero $a \in \ff_1$, then $b =0$.
    \item[(3)] The truncated Lie algebra $\ff/(\bigoplus_{k=2}^{\infty}\ff_{-k}) = \ff_{-1} \oplus \ff_0 \oplus \ff_1$ is isomorphic to $\fs\fl_2$.  \end{itemize} \end{lemma}

We have the following partial generalization of Theorem \ref{t.Yamaguchi}.

\begin{lemma}\label{l.Ya}
Let $\fl = \fl_{-1} \oplus \fl_0 \oplus \fl_1$ be the graded Lie algebra associated with a subadjoint variety $Z = L/P \subset \BP \bV$ in Definition \ref{d.fh}. Then  $\fl_{-1}$ is the first prolongation of $\fl_0 \oplus \fl_1$.
\end{lemma}

\begin{proof}
 This is immediate from Theorem \ref{t.Yamaguchi} in the case (i) of Proposition \ref{p.list}.

In the case of (ii-1) of Proposition \ref{p.list},  we can write $\fl = \fq' \oplus \fq''$ where $\fq' \cong  \fq'' \cong \fsl_2$ has the grading coming from $\BP^1$ and derive the statement of the lemma from Lemma \ref{l.sumprolong}  and Lemma \ref{l.sl2}.

In the case of (ii-2) of Proposition \ref{p.list}, we can write $\fl = (\fq' \oplus \fq'' ) \oplus \fq'''$ where $\fq' \cong \fq'' \cong \fq''' \cong \fsl_2(\C)$ with the grading coming from $\BP^1$ and deduce the statement of the lemma  by applying  Lemma \ref{l.sumprolong}  and Lemma \ref{l.sl2} twice.

In the case of  (ii-3) of Proposition \ref{p.list}, we can write $\fl = \fq' \oplus \fq''$ where $\fq' \cong  \fsl_2$ has the grading coming from $\BP^1$ and $\fq'' \cong \fs\fo_{m+4}$  has the grading coming from the quadric hypersurface $\Q^m, m \geq 3$. Thus the statement of the lemma follows from Lemma \ref{l.sumprolong} combined with Theorem \ref{t.Yamaguchi} for $\Q^m, m \geq 3$.
\end{proof}

We prove Theorem \ref{t.prolong} in two parts: Proposition \ref{p.prolong1} for $\fp_{-1} \cong \fg_{-1}$ and Proposition \ref{p.prolong2} for $\fp_{-2} =0$. Each of the two propositions is preceded by some lemmata.

\begin{lemma}\label{l.prolong11}
In the setting of Theorem \ref{t.prolong},
for any $\varphi \in \fp_{-1} \subset \Hom(\fg_+, \fg)_{-1}$ and $a, b \in \fl_1$, the following holds in $\fg$.
\begin{itemize}
\item[(i)] $[\varphi(a), b] = [\varphi(b), a]$;
\item[(ii)] $[ \varphi(a \cdot v_0), b \cdot v_0] = [ \varphi(b \cdot v_0), a \cdot v_0]$; and
\item[(iii)] $[\varphi(a), b \cdot v_0] + [\varphi(a \cdot v_0), b] = [\varphi(b), a \cdot v_0] + [\varphi(b \cdot v_0), a].$ \end{itemize} \end{lemma}

    \begin{proof} Since $[a, b] = [a\cdot v_0, b \cdot v_0] =0$ for all $a, b \in \fl_1$,
    (i) (resp. (ii)) follows from (\ref{e.prolong}) by putting $x=a, y=b$ (resp. $x = a \cdot v_0, y = b \cdot v_0$).
   (iii) follows from (i), (ii), (\ref{e.prolong}) and (\ref{e.II})    by putting $x = a + t a \cdot v_o, y = b+ t b \cdot v_o.$
    \end{proof}

    \begin{lemma}\label{l.prolong12}
    In Lemma \ref{l.prolong11}, let $\lambda, \mu \in \fg_1^{\vee}$ and $\xi \in \Hom(\fg_1, \fl_0)$ such that
    \begin{equation}\label{e.l12} \varphi(x) = \lambda(x) v_0 + \mu(x) {\rm Id}_{\bV} + \xi(x) \ \in \ \bV_o \oplus \C {\rm Id}_{\bV} \oplus \fl_0 = \fg_0 \end{equation}  for all $x \in \fg_1$.  Then for any $a, b \in \fl_1$, \begin{itemize}
    \item[(1)] $\lambda (a) b \cdot v_0 = \lambda(b) a \cdot v_0 $;
    \item[(2)] $[\xi(a), b] = [\xi(b), a]$;
    \item[(3)] $[\xi(a \cdot v_0), b ] = [\xi(b \cdot v_0), a]$;
    \item[(4)] $\mu(a \cdot v_0) b\cdot v_0 + [\xi(a \cdot v_0), b\cdot v_0] = \mu(b \cdot v_0) a \cdot v_0 + [\xi(b \cdot v_0), a \cdot v_0]$; and
    \item[(5)] $(\mu(a) - \lambda(a \cdot v_0)) b\cdot v_0 + [ \xi(a), b \cdot v_0] = (\mu(b) - \lambda(b \cdot v_0)) a \cdot v_0 + [\xi(b), a \cdot v_0].$
        \end{itemize}
        It follows that $\lambda (\fl_1) =0$ and there exist unique elements $\theta, \zeta \in \fl_{-1}$ such that $$\xi(a) = [\theta, a] \mbox{ and } \xi(a \cdot v_0) = [\zeta, a]$$ for all $a \in \fl_1.$ \end{lemma}

        \begin{proof}
        (1)-(5) follow from (i)-(iii) of Lemma \ref{l.prolong11} by substituting (\ref{e.l12}) and separating $\bV$-components and $\fl$-components. Then
        (1) implies $\lambda (\fl_1) = 0$ because $\dim \fl_1 >1.$  It remains to show the existence and the uniqueness of $\theta, \zeta$ in the last sentence of the lemma. But this follows from Lemma \ref{l.Ya} because
        the homomorphism $\fl_1 \to \fl_0$ defined by $a \mapsto \xi(a)$ (resp. $a \mapsto \xi(a \cdot v_0)$)
        is an element of the first prolongation of $\fl_1 \oplus \fl_0$ by (2) (resp. by (3)).
 \end{proof}

\begin{proposition}\label{p.prolong1}
In the setting of Theorem \ref{t.prolong},  we have the following.
\begin{itemize}
\item[(1)] The Lie bracket of $\fg$ induces an inclusion $\fl_{-1} \subset \fp_{-1}.$
\item[(2)] The homomorphism $\fp_{-1} \to \fl_{-1}$ given by $\varphi \mapsto \theta$ in the notation of Lemma \ref{l.prolong12} is injective. \end{itemize}
    In particular, we have an isomorphism $\fl_{-1} \cong \fp_{-1}$. \end{proposition}

\begin{proof}
The Lie bracket induces a homomorphism $$\fl_{-1} \to \Hom(\fg_+, \fg_0 \oplus \fg_+)_{-1}$$ whose image is contained in $\fp_{-1}$. This is injective because (e.g. by the relation of the grading with the root space decomposition in Section 3.3 of \cite{Ya}) $$\{ a \in \fl_{-1} \mid  [a, \fl_1] =0\} =0,$$ which proves (1).

To prove (2), let $\varphi$ be an element of  $\fp_{-1}$ satisfying $\theta =0$. We want to show that $\varphi=0$.
 Putting $\xi(a) = [\theta, a] = 0$ in  Lemma \ref{l.prolong12} (5), we have $$
 \mu(a) = \lambda( a \cdot v_0) \mbox{ for all } a \in \fl_1.$$
 Using the notation of Definition \ref{d.fh}, put  $$ [\xi(a \cdot v_0), b \cdot v_0] = {\bf c}([\zeta, a]) b \cdot v_0+[[\zeta, a], b]\cdot v_0  \mbox{ and } $$ $$  [\xi(b \cdot v_0), a \cdot v_0] = {\bf c}([\zeta, b]) a \cdot v_0+[[\zeta, b], a]\cdot v_0$$  into Lemma \ref{l.prolong12} (4) to obtain $$( \mu(a \cdot v_0) + {\bf c}([\zeta, a]) ) b\cdot v_0 = (\mu(b \cdot v_0) + {\bf c}([\zeta, b]) a \cdot v_0$$ for all $a, b \in \fl_1$. Since $\dim \fl_1 >1$, this implies $\mu(a \cdot v_0) = - {\bf c}([\zeta, a])$ for all $a \in \fl_1$. In summary, we have 
 \begin{equation}\label{e.pr1} \varphi(a) = \lambda(a \cdot v_0) {\rm Id}_{\bV} \end{equation} \begin{equation}\label{e.pr2} \varphi( a \cdot v_0) = \lambda(a\cdot v_0)v_0 - {\bf c}([\zeta, a]) {\rm Id}_{\bV} + [\zeta, a] \end{equation} 
 for all $a \in \fl_1.$
Consequently, we obtain  
 \begin{equation}\label{e.pr3} \varphi(a^2 \cdot v_0) = [\varphi(a), a \cdot v_0] + [ a, \varphi(a \cdot v_0)] = 2 \lambda(a\cdot v_0) a \cdot v_0 + [a, [\zeta, a]], \end{equation} \begin{equation}\label{e.pr4} \varphi( a^3 \cdot v_0) = [\varphi(a), a^2 \cdot v_0] + [a, \varphi(a^2 \cdot v_0)] = 3 \lambda(a \cdot v_0) a^2 \cdot v_0, \end{equation} 
 where we used $[a, [a, [\zeta, a]]] \in [\fl_1, \fl_1] =0$ to derive (\ref{e.pr4}).

 Consider \begin{eqnarray*} \widehat{{\rm Bs}}(\II) & = & \{ a \in \fl_1 \mid \II(a, a) =0\} \mbox{ and } \\  \widehat{{\rm Bs}}(\III) &=& \{ a \in \fl_1 \mid \III(a, a,a) =0\}. \end{eqnarray*}  If $a \in \widehat{{\rm Bs}}(\II)$, then $a^2 \cdot v_0 = {\rm II}(a,a) \cdot v_0 =0$. Thus \begin{equation}\label{e.pr5} \lambda(a\cdot v_0) = 0 = [[\zeta,a],a] \mbox{ for all } a \in \widehat{{\rm Bs}}({\rm II}) \subset \fl_1.\end{equation}

 If $Z \subset \BP V$ is not the case of (ii-1) of Proposition \ref{p.list}, then $ \widehat{{\rm Bs}}({\rm II})  = \widehat{\sB}$ in Proposition \ref{p.FFsubadjoint} (5) is linearly nondegenerate in $\fl_1$. Thus (\ref{e.pr5}) implies that $\lambda (\bV_1) =0$. Moreover,  Lemma \ref{l.xvv} and (\ref{e.pr5}) show $\zeta =0$, proving $\varphi =0$.

 It remains to handle  the case (ii-1) of Proposition \ref{p.list} with $Z  \cong \BP^1 \times \Q^1$. In this case, the variety $\widehat{{\rm Bs}}({\rm II})$ is linearly degenerate in $\fl_1$, corresponding to the $\BP^1$-factor in $Z $, while $\widehat{{\rm Bs}}({\rm III})$
 is  nondegenerate  in $\fl_1$.
 For $a \in \widehat{{\rm Bs}}({\rm III}) \setminus \widehat{{\rm Bs}}({\rm II})$, (\ref{e.pr4}) implies that $\lambda(a\cdot v_0) =0$. It follows that $\lambda(\bV_1) =0$. To prove $\zeta=0$, write $\fl= \fl'\oplus \fl''$ with the grading $\fl'= \fl'_{-1} \oplus \fl'_0 \oplus \fl'_1$ given by $\BP^1$ and $\fl'' = \fl''_{-1} \oplus \fl''_0 \oplus \fl''_{-1}$ given by $\Q^1 \cong \BP^1$. We have $\zeta' \in \fl'_{-1}$ and $\zeta'' \in \fl''_{-1}$ satisfying $\zeta= \zeta' + \zeta'' \in \fl_{-1}.$  From the second fundamental form of $\BP^1 \times \Q^1 \subset \BP^5$, we see that
 ${\rm II}(\fl'_1, \fl'_1) =0$ and
 \begin{equation}\label{e.pr6}
 {\rm II}(a,b) \neq 0 \mbox{ for any nonzero } a, b \in \fl''_1. \end{equation}
 Since $[\zeta'', \fl'_1] =0$ and $\widehat{{\rm Bs}}({\rm II}) = \fl'_1$, we know from (\ref{e.pr3}) that $$[[\zeta', a], a] =0 \mbox{ for any } a \in \fl'_1.$$
 From the property of $\fl' \cong \fs\fl _2$ (or applying Lemma \ref{l.xvv} to $\fl' = \fl'_{-1} \oplus \fl'_0 \oplus \fl'_{-1}$), this implies $\zeta'=0$ and $\zeta = \zeta''.$ In summary, we have for any $a \in \fl_1$,
 $$\varphi(a) = \lambda(a \cdot v_0) {\rm Id}_{\bV} = 0,$$
 $$\varphi(a \cdot v_0) = - {\bf c}([\zeta'',a]) {\rm Id}_{\bV} + [\zeta'', a] \in \C {\rm Id}_{\bV} \oplus \fl''_0,$$
 $$\varphi(a^2 \cdot v_0) = [a, [\zeta'', a]] \in \fl''_1, $$
 $$\varphi(a^3 \cdot v_0) = 3 \lambda(a \cdot v_0) a^2 \cdot v_0 =0.$$
  Using $[\zeta'', a] \cdot v_0 = {\bf c}([\zeta'', a]) v_0$ for any $a \in \fl_1$, we have
 \begin{equation}\label{e.pr7} [\varphi(a \cdot v_0), a^2 \cdot v_0] = 2 [[\zeta'', a], a] a \cdot v_0 = 2 {\rm II}([[\zeta'',a],a], a) \cdot v_0, \end{equation}
\begin{equation}\label{e.pr8} [a \cdot v_0, \varphi(a^2\cdot v_0)] = [[\zeta'', a], a] \cdot (a \cdot v_0) = {\rm II}([[\zeta'', a], a], a) \cdot v_0. \end{equation}
 From $[\bV_1, \bV_2] =0$, we have $$0 = \varphi([a\cdot v_0, a^2\cdot v_0])
 = [\varphi(a \cdot v_0), a^2 \cdot v_0] + [a \cdot v_0, \varphi(a^2 \cdot v_0)].$$
Thus (\ref{e.pr7}) and (\ref{e.pr8}) give
${\rm II}( [[\zeta'', a],a], a) =0$ for all $a \in\fl_1$. It follows from (\ref{e.pr6})
that $[[\zeta'', a],a] =0$ for all $a \in \fl''_1$. Thus by the property of $\fs\fl_2$
(or applying Lemma \ref{l.xvv} to $\fl''$), we see that $\zeta''=0$, completing the proof that $\varphi =0$. \end{proof}

Turning to the second prolongation, we have the following lemmata.

\begin{lemma}\label{l.prolong2}
In the setting of Theorem \ref{t.prolong}, for an element $\varphi \in \fp_{-2} \subset  \Hom(\fg_{+}, \fg)_{-2}$, let  $\phi,
\psi \in \Hom(\fl_1, \fl_{-1})$ be the homomorphisms satisfying $\varphi(a) = \phi(a) \in \fl_{-1}$ and $\varphi(a\cdot v_o) = \psi(a)\in\fl_{-1}$ for all $a \in \fl_1.$ Then for any $a, b \in \fl_1,$
\begin{itemize}
\item[(1)] $[\phi(a), b]  =  [\phi(b), a]$;
\item[(2)] $[\psi(a), b\cdot v_0] = [\psi(b), a\cdot v_0]$; and
\item[(3)] $[\phi(a), b\cdot v_0] + [\psi(a), b] = [\phi(b), a\cdot v_0] + [\psi(b), a].$
\end{itemize}
In particular, the two elements $\phi, \psi \in \Hom(\fl_1, \fl)_{-2}$ lie in the second prolongation of $\fl_0 \oplus \fl_1$.
\end{lemma}

\begin{proof}
Put $x = a+ t a\cdot v_0,  y = b + t b \cdot v_0$ in (\ref{e.prolong}) and apply (\ref{e.II})
to obtain $$[\varphi(a + ta \cdot v_0), b + t b \cdot v_0] = [\varphi(b + t b \cdot v_0), a + t a \cdot v_0]$$
for all $t \in \C$.
(1), (2) and (3) are obtained by comparing the  coefficients of the above polynomial equation in $t$.
(1) implies that $\phi$ is in the second prolongation of $\fl_1 \oplus \fl_0$.
Considering the $\fl$-components in (3), we have $[\psi(a), b] = [\psi(b), a]$ for all $a, b \in \fl_1$. Thus $\psi$ is also in the second prolongation of $\fl_1 \oplus \fl_0$.
\end{proof}

\begin{lemma}\label{l.sl2W}
Let $\ff_{-1} \oplus \ff_0 \oplus \ff_1$ be the gradation on $\fs\fl_2$ associated with $\BP^1$. Let $W$ be the adjoint representation of $\fs\fl_2$ and $w_0 \in W$ be the lowest weight vector.  For any nonzero element $\phi \in \ff_{-2}$ of the second prolongation of $\ff_0 \oplus \ff_1$ and any nonzero element $a \in \ff_1$, we have
$$[[ \phi(a), a], a] \cdot w_0 = -a \cdot ([\phi(a), a] \cdot w_0)  \neq 0.$$
\end{lemma}

\begin{proof}
In the notation of Lemma \ref{l.sl2},
we may assume $$w_0 = t^2\frac{\p}{\p t}, \ a = \frac{\p}{\p t}, \ \phi = [t^3 \frac{\p}{\p t}, \cdot].$$  Then
$$[\phi(a), a] = [ [ t^3 \frac{\p}{\p t}, \frac{\p}{\p t}], \frac{\p}{\p t}] = 6 t \frac{\p}{\p t},$$ $$
[[\phi(a), a], a] = [6 t \frac{\p}{\p t}, \frac{\p}{\p t}] = -6  \frac{\p}{\p t}.$$ Thus
$$[[\phi(a), a], a] \cdot w_0 = -6  [\frac{\p}{\p t}, t^2 \frac{\p}{\p t}] = -12t \frac{\p}{\p t},$$
$$a \cdot ([\phi(a), a] \cdot w_0) = [\frac{\p}{\p t}, [6t \frac{\p}{\p t}, t^2 \frac{\p}{\p t}]] = [\frac{\p}{\p t}, 6 t^2 \frac{\p}{\p t}] = 12t  \frac{\p}{\p t} $$
which proves the lemma. \end{proof}

\begin{proposition}\label{p.prolong2}
In the setting of Theorem \ref{t.prolong}, we have $\fp_{-k} =0$ for $k \geq 2$.
\end{proposition}

\begin{proof}
It suffices to prove $\fp_{-2} =0$.
Using the terminology of Lemma \ref{l.prolong2}, we need to show that $\phi = \psi =0$.
When $\fl= \fl_1 \oplus \fl_0 \oplus \fl_{-1}$ is associated with the case (i) in Proposition \ref{p.list}, we have $\phi= \psi =0$ by
Theorem \ref{t.Yamaguchi}.

Thus we may assume that $\fl$ is associated with the case (ii)  in Proposition \ref{p.list} corresponding to the Segre embedding $\BP^1 \times \Q^{m-2} \subset \BP^{2m-1}$ for some $m \geq 3$. We can write $\fl= \fl' \oplus \fl''$ such that
$\fl' = \fl'_{-1} \oplus \fl'_0 \oplus \fl'_1$ is the grading coming from $\BP^1$ and $\fl'' = \fl''_{-1} \oplus \fl''_0 \oplus \fl''_1$ is the grading coming from $\Q^{m-2}$. Let $\fl'_{-k}$ (resp. $\fl''_{-k}$) be the $k$-th prolongation of $\fl'_0 \oplus \fl'_1$ (resp. $\fl''_0 \oplus \fl''_1$).
Both $\phi$ and $\psi$ are in the second prolongation of $\fl_0 \oplus \fl_1$, which is isomorphic to $\fl'_{-2} \oplus \fl''_{-2}$ by Lemma \ref{l.sumprolong}. So there exist $\phi', \psi' \in \fl'_{-2}$ and $\phi'', \psi'' \in \fl''_{-2}$ such that
$\phi = \phi' + \phi''$ and $\psi= \psi' + \psi''$.

We claim that $\phi'= \psi'=0$. Since the $\BP^1$-factor of $Z = \BP^1 \times \Q^{m-2}$ is linearly embedded in $\BP \bV$, the subspace $\fl'_1$ of $\fl_1$ is in the base locus $\widehat{{\rm Bs}}({\rm II})$ of the second fundamental form of $Z$. This means
$[a, a \cdot v_0] =0$ for any $a \in \fl'_1$. Hence \begin{equation}\label{e.l2}
[\varphi(a), a \cdot v_0] = [ \varphi(a \cdot v_0), a] \mbox{ for any } a\in\fl'_1. \end{equation}
The left-hand side of (\ref{e.l2}) is \begin{eqnarray*}
[\varphi(a), a \cdot v_0] &=& [\phi'(a), a \cdot v_0]  \\
&=& [[\phi'( a), a], v_0] + [a, [\phi'(a), v_0]] \\
&=& [\phi'(a), a] \cdot v_0 \ \in \bV_0 \end{eqnarray*}
where the last equality used $[\phi'(a), v_0] =0$  coming from the fact that $v_0$ is the lowest weight vector and $\phi'(a) \in \fl'_{-1}$. The right-hand side of (\ref{e.l2}) is
$$[\varphi(a \cdot v_0), a] = [\psi'(a), a]  \in \fl'_0.$$
Since the left-hand side of (\ref{e.l2}) is in $\bV_0$ and the right-hand side is in $\fl'_0$, we see that both sides are zero.
Thus $[\phi'(a), a] = [\psi'(a), a] =0$ for any $a \in \fl'_1$, which implies  $\phi'= \psi' =0$ by Lemma \ref{l.sl2} (2).

It remains to check $\phi'' = \psi'' =0$. This is immediate from Theorem \ref{t.Yamaguchi} when $Z \cong \BP^1 \times \Q^{m-2}$ with $m\geq 5$.
When $m=4,$ our proof for $\phi'= \psi'=0$ can be applied to show $\phi''=\psi''=0$  by exchanging the $\BP^1$-factors in $Z \cong \BP^1 \times \BP^1 \times \BP^1.$
Thus we may assume that $m=3, \fl'' \cong \fs\fl_2$ and $Z \cong \BP^1 \times \Q^1 \subset \BP^5.$
Then the subspace $\fl''_1$ of $\fl_1$ is in $\widehat{{\rm Bs}}({\rm III})$ and we have
$[a, a^2 \cdot v_0] =0$ for any $a\in \fl''_1$. This implies
\begin{equation}\label{e.l3}
[\varphi(a), a^2 \cdot v_0] = [\varphi(a^2 \cdot v_0), a]. \end{equation}
Since $[\phi''(a), v_0 ] =0$ for the lowest weight vector $v_0$ and any $a\in \fl''_1$, we have \begin{eqnarray*} \varphi(a^2 \cdot v_0) &=& [\varphi(a), a\cdot v_0] + [ a, \varphi(a \cdot v_0)]\\ &=& [\phi''(a), a \cdot v_0] + [a, \psi''(a)] \\
&=& [\phi''(a), a] \cdot v_0 + [a, \psi''(a)]. \end{eqnarray*}
Thus the right-hand side of (\ref{e.l3}) is
\begin{equation}\label{e.l4} [[[\phi''(a), a] \cdot v_0] , a] - [ [\psi''(a), a], a].
\end{equation}
On the other hand, the left-hand side of (\ref{e.l3}) is
\begin{eqnarray*} [\phi''(a), [a, [a, v_0]] ] &=& [[\phi''(a), a], [a, v_0]] + [a, [ \phi''(a), [a, v_0]]] \\ &=&
  [[[\phi''(a), a], a], v_0] + [a, [[\phi''(a), a],  v_0]]] \\ & & + [a,[ [\phi''(a),a], v_0] ] \\
 &=&  [[[\phi''(a), a], a], v_0] + 2 [a, [[\phi''(a), a],  v_0]]], \end{eqnarray*}
 where we used   $[\phi''(a), v_0] =0$.
Equating it with (\ref{e.l4}), we obtain
$$3 [a, [[\phi''(a), a],  v_0]]] + [[[\phi''(a), a], a], v_0] = - [[ \psi''(a), a], a].$$ Since the left-hand side is in $\bV_1$ while the right-hand side is in $\fl_1$, both sides vanish. Thus $ \psi''(a) =0$  from Lemma \ref{l.sl2} (2) and $\phi''=0$ from Lemma \ref{l.sl2W},  because $\Q^1 \subset \BP^2$ is the generalized minuscule variety associated to the adjoint representation of $\fs\fl_2$.
\end{proof}

\section{Cohomological bundles on subadjoint varieties}\label{s.Spencer}
In this section and the next section, we examine the vanishing condition of Theorem \ref{t.HL} for the graded Lie algebra $\fg$ in Definition \ref{d.fh}.

\begin{notation}\label{n.EZ}
Let $Z = L/P$ be a generalized minuscule variety in Notation \ref{n.hss}. Given an $L_0$-module $V$, we regard it as a $P$-module via the projection $P \to P/L_{-1} = L_0$ and denote by $\sE_Z(V)$ the vector bundle on $Z$ associated to the $P$-principal bundle $L \to L/P = Z.$ \end{notation}

The goal of this section is to describe a distinguished $L_0$-submodule
$R_k \subset Q^{k}(\fg)$ for the Lie algebra $\fg$ of Definition \ref{d.fh} such that  $$ H^0(Z, \sE_Z(Q^k(\fg)))= H^0(Z, \sE_Z(R_k)).$$ We start the discussion by recalling some standard results on homogeneous vector bundles, the formulation of which needs the following terminology.

\begin{notation}\label{n.bc} Let $\fl$ be a semisimple Lie algebra. Fix a Cartan subalgebra  $\ft$ and a set of simple roots $\alpha_1, \ldots, \alpha_{\ell}$ with respect to $\ft$.
\begin{itemize} \item[(1)] For  a weight $\omega \in \ft^{\vee}$, we can write
$$\omega=\sum_{k=1}^{\ell} \bc^{\alpha_k}(\omega) \alpha_k$$ where $\bc^{\alpha_k}(\omega)$ is a rational number  determined by this expression.
 \item[(2)] For $I$ as in Notation \ref{n.hss} and a weight $\omega \in \ft^{\vee},$   define $$\bc^I(\omega)=\sum_{i\in I}\bc^{\alpha_i}(\omega)$$ which is a rational number.
\item[(3)] Given any $\ft$-module $V$, we denote by ${\rm Weights}(V) \subset \ft^{\vee}$  the set of weights of $V$ and define
\begin{eqnarray*}
\bc^I(V):=\{\bc^I(\omega)\mid \omega\in {\rm Weights} (V)\}
\end{eqnarray*}
which is a set of rational numbers. \end{itemize}
\end{notation}

\begin{lemma}\label{l.BWB}
In Notation \ref{n.EZ} and Notation \ref{n.bc}, if  $H^0(Z, \sE_Z(V)) \neq 0$, then $\bc^I(V)$ contains a nonnegative rational number.
\end{lemma}

\begin{proof} Let  $V=\bigoplus_{j=1}^{k}V_j$ be  a decomposition into  irreducible $L_0$-modules.  Then $$H^0(Z, \sE_Z(V))=\bigoplus_{j=1}^k H^0(Z, \sE_Z(V_j)).$$
If $H^0(Z, \sE_Z(V))\neq 0$, then $H^0(Z, \sE_Z(V_{j_0}))\neq 0$ for some $1 \leq j_0 \leq k$.  By Borel-Weil-Bott theorem (e.g. Section 4.3 of \cite{Ak}), the highest weight $\omega$ of $V_{j_0}$ is dominant. Thus $\bc^{\alpha_i}(\omega)$ is a nonnegative rational number for each $i\in I$, and so is $\bc^I(\omega)$.
\end{proof}

\begin{lemma}\label{l.genminuscule}
Let  $Z  \subset\BP V(\sum_{i \in I} d_i \omega_i)^{\vee}$ be a generalized minuscule variety  for some positive integers $d_i, i \in I$.  Let $$v_0 \in V(- \bw_0(\sum_{i\in I} d_i \omega_i)) = V(\sum_{i \in I} d_i \omega_i)^{\vee}$$ be a lowest weight vector and let $V_j, 0 \leq j \leq r,$ be as in Definition \ref{d.minuscule}.  Then
\begin{itemize} \item[(1)] $\bc^I(\fl_j) = \{ j \}$ for all $j = -1, 0, 1$;
\item[(2)] $\bc^I(V_j) = j - \bc^I(\sum_{i\in I} d_i \omega_i)$ for all $0 \leq j \leq r$.
    \end{itemize}
\end{lemma}

\begin{proof}
(1) follows from the fact that $\fl_j =\bigoplus_{i\in I}\fl_j^i $ for $j= -1, 0, 1,$ and $\bc^{\alpha_k}(\fl_j^i)=\{\delta_{ik}j\}$, where $i, k\in I$ and $\delta_{ik}$ is the Kronecker  symbol.
(2) follows from  $$V_j=\Sym^j\fl _1\cdot V_0, \ \  {\rm Weights} (V_0)= \{ -\sum_{i\in I} d_i \omega_i\},$$ and $\bc^I(\fl_1)=\{ 1\}$ from (1).
\end{proof}

The following  lemma  can be checked easily from Table 1 of \cite{Hum}.

        \begin{lemma}\label{l.bcsubadjoint}  When $\sO(\omega_*)$ is the embedding line bundle of the subadjoint variety $Z \subset \BP \bV$ in Definition \ref{d.fh}, the number $\bc^I(\omega_*)$ is $\frac{3}{2}$.
\end{lemma}

To state the main result of this section, Proposition \ref{p.R_k}, we need  the following three lemmata.

\begin{lemma}\label{l.decompo}
In Definition \ref{d.fh},  set $$\wg_{-1} = \fl_{-1}, \ \wg_0 = \C {\rm Id}_{\bV} \oplus \fl_0, \ \wg_1 = \fl_1$$ and $\wg = \wg_{-1} \oplus \wg_0 \oplus \wg_1$ such that $\fg = \wg \ltimes \bV.$ For each negative integer $k$,  we have the following direct sum decomposition of $L_0$-modules
\begin{eqnarray*}
\Hom(\wedge^2 \fg_+, \fg)_k & = & \Hom( \wedge^2 \wg_1, \bV_{k+2}) \oplus \Hom(\wedge^2 \wg_1, \wg_{k+2})  \\ & &
\oplus \bigoplus_{i \geq 1} \Hom(\wg_1 \otimes \bV_i, \wg_{k+i+1}) \\ & &
\oplus \bigoplus_{i \geq 1} \Hom(\wg_1 \otimes \bV_i, \bV_{k+i+1}) \\ & &
\oplus \bigoplus_{j \geq i \geq 1} \Hom( \bV_i \wedge \bV_j, \wg_{k + i +j}) \\ & &
\oplus \bigoplus_{j \geq i \geq 1} \Hom(\bV_i \wedge \bV_j, \bV_{k + i + j}). \end{eqnarray*}

\end{lemma}

\begin{lemma}\label{l.partial}
For the universal prolongation $\fg$ of Theorem \ref{t.prolong}, the image of the homomorphism $\partial: C^{-1,1}(\fg) \to C^{-1,2}(\fg)$ in Definition \ref{d.cohomology} contains an $L_0$-submodule  isomorphic to $\Hom(\wedge^2 \bV_2, \bV_3)$.
\end{lemma}

\begin{proof}
Consider the  two homomorphisms defined by the restrictions of $\partial: C^{-1,1}(\fg) \to C^{-1,2}(\fg),$
\begin{eqnarray*}
&& \partial': \Hom(\bV_2, \fl_1) \to \Hom(\wedge^2 \bV_2, \bV_3), \\
&& \partial'': \Hom(\bV_2, \fl_1) \to \Hom(\bV_1\wedge \bV_2, \bV_2).
\end{eqnarray*}
We claim that $\partial'$ is surjective and $\partial''$ is injective.
Since both  are homomorphisms of $L_0$-modules and $L_0$ is reductive, the claim implies the existence of an injective homomorphism of $L_0$-modules $$\Hom(\wedge^2 \bV_2, \bV_3)\to\partial(\Hom(\bV_2, \fl_1)),$$ which proves the lemma.

To show the claim, note that the Lie bracket $\bV_2 \otimes \fl_1 \to \bV_3 \cong \C $ is a perfect pairing from  Proposition \ref{p.FFsubadjoint} (4) and Lemma \ref{l.II}. By this perfect pairing, the  homomorphism $\partial'$ is just the anti-symmetrization $\bV_2^{\vee} \otimes \bV_2^{\vee} \to \wedge^2 \bV_2^{\vee}$, which is surjective.
Define $L_0$-modules $$S_1:=\{f\in\Hom(\bV_2, \fl_1)\mid\partial f(\bV_1\wedge\bV_2)=0\} \mbox{ and } $$ $$S_2:=\{f(v)\in\fl_1\mid f\in S_1, v\in \bV_2\}.$$ For any $f\in S_1$, $v_1\in\bV_1$ and $v_2\in\bV_2$, we have $$f(v_2)\cdot v_1=-\partial f(v_1, v_2)=0,$$ implying $[S_2, \bV_1]=0$. Since $[\fl_1, \bV_1]=\bV_2$, the $L_0$-submodule $S_2$ of the irreducible $L_0$-module $\fl_1$ must be the zero module. It follows that $f(\bV_2)=0$ for all $f\in S_1.$ This shows that $S_1=0$ and  the homomorphism $\partial''$ is  injective.
\end{proof}

The following lemma is  immediate since $L_0$ is reductive.

\begin{lemma}\label{l.iota} In the setting of Lemma \ref{l.partial},
the surjective $L_0$-module homomorphism
$$C^{k,2}(\fg) = \Hom(\wedge^2 \fg_+, \fg)_k \to Q^{k}(\fg) = C^{k,2}(\fg)/\partial (C^{k,1}(\fg)) $$
admits an $L_0$-module splitting $Q^{k}(\fg) \subset \Hom(\wedge^2 \fg_+, \fg)_k $ for each negative integer $k$. Furthermore, when $k=-1$, Lemma \ref{l.partial} gives an injective homomorphism of $L_0$-modules
$$Q^{-1}(\fg) \to  \Hom(\wedge^2 \fg_+, \fg)_{-1}/\Hom(\wedge^2 \bV_2, \bV_3).$$
They  induce    injective homomorphisms \begin{eqnarray*} \lefteqn{ \iota_{-1}:
H^0(Z, \sE_Z(Q^k(\fg))) \to} \\ & &  H^0(Z, \sE_Z(\Hom(\wedge^2 \fg_+, \fg)_{-1}/\Hom(\wedge^2 \bV_2, \bV_3))) \end{eqnarray*}  and   $$ \iota_k:
H^0(Z, \sE_Z(Q^k(\fg))) \to H^0(Z, \sE_Z(\Hom(\wedge^2 \fg_+, \fg)_k))$$ for $k \leq -2$.
\end{lemma}

\begin{proposition}\label{p.R_k}
In the setting of Lemma \ref{l.iota}, define \begin{eqnarray*}
R_{-1} & := & \Hom(\bV_1 \wedge \bV_+, \bV)_{-1} \oplus \Hom(\wedge^2\bV_1, \wg_1) \oplus \Hom(\fl_1 \wedge \bV_1, \wg_1); \\
R_k & := &  \Hom ( \wedge^2 \bV_+, \wg)_k \mbox{ for } k = -2, -3; \mbox{ and } \\
R_k & := & 0 \mbox{ for } k \leq -4. \end{eqnarray*} Then the image of $\iota_k$ in Lemma \ref{l.iota} is contained in $H^0(Z, \sE_Z(R_k)).$ \end{proposition}

\begin{proof}
We need to show that $H^0(Z, \sE_Z(V)) =0$ for any irreducible $L_0$-submodule  $V$ of
$$\Hom(\wedge^2 \fg_+, \fg)_{-1}/\Hom(\wedge^2 \bV_2, \bV_3) \mbox{ or }\Hom(\wedge^2 \fg_+, \fg)_{k \leq -2}$$  satisfying  $V \cap R_k =0.$
By Lemma \ref{l.BWB}, it suffices to show that $\bc^I(V)$ does not contain nonnegative rational numbers unless $V \cap R_k \neq 0.$ We check this by applying the decomposition in  Lemma \ref{l.decompo}.
Using Lemma \ref{l.genminuscule} and Lemma \ref{l.bcsubadjoint}, we have
\begin{eqnarray*}
\bc^I(\Hom(\wedge^2 \wg_1, \bV_{k+2})) &=& \{ k- \bc^I(\omega_*)\} = \{ k- \frac{3}{2}\}, \\
\bc^I(\Hom(\wedge^2 \wg_1, \wg_{k+2})) & = & \{ k \}, \\
\bc^I(\Hom(\wg_1 \otimes \bV_+, \bV)_k) &=& \{ k \}, \\
\bc^I(\Hom(\wg_1 \otimes \bV_+, \wg)_k) & =& \{  k + \bc^I(\omega_*)\} = \{ k + \frac{3}{2}\}, \\
\bc^I(\Hom(\wedge^2 \bV_+, \bV)_k) &=& \{ k + \bc^I(\omega_*) \} = \{ k + \frac{3}{2}\}, \\
\bc^I(\Hom(\wedge^2 \bV_+, \wg)_k) & = & \{ k + 2 \bc^I(\omega^*)\} = \{ k + 3\}.\end{eqnarray*}
It is easy to see that the sets of rational numbers on the right hand side can contain nonnegative rational numbers only when the $L_0$-modules on the left hand side have positive-dimensional intersection with  $R_{-1} + \Hom(\wedge^2 \bV_2, \bV_3)$, $R_{-2}$, or $R_{-3}$.
\end{proof}

\section{Cohomological bundles associated with $\bG_0$-principal bundles }\label{s.M}

\begin{definition}\label{d.symbol}
Let $D \subset TY$ be a distribution on a manifold $Y$, i.e., a vector subbundle of the tangent bundle.
Let $\sO(D)$ be the locally free subsheaf of the sheaf $\sO(TY)$ of vector fields on $Y$ and denote by $\partial \sO(D) = \partial^1 \sO(D)$ the saturated subsheaf of $\sO(TY)$ generated by $ \sO(D) + [\sO(D), \sO(D)]$. Define $\partial^{i+1} \sO(D)$ inductively  as the saturated subsheaf of $\sO(TY)$ generated by $\partial^i \sO(D) + [ \partial^i \sO(D), \sO(D)].$ We say that $D$ is a {\em regular distribution} if there exists a vector subbundle $\partial^i D \subset TY$ for each $i \geq 1$ such that $\partial^i \sO(D)$ is the locally free sheaf associated to $\partial^i D$. By convention, we put $\partial^0 D = D$. For a regular distribution $D$, consider the graded vector space defined by
$${\rm Symb}_y D := \bigoplus_{i \geq 0} (\partial^{i} D)_y / (\partial^{i-1} D)_y,$$ for each $y \in Y$. This  defines a graded vector bundle ${\rm Symb}D$ on $Y$. Note that a regular distribution $D$
determines a filtration $F^{\bullet}$ on $Y$ by $F^i = \partial^{i-1} D$ and ${\rm Symb}_y D$ is isomorphic to the symbol algebra of the filtration, which we call
the {\em symbol algebra} of $D$ at $y$. \end{definition}

In this section, we work in the following setting.

\begin{setup}\label{setup}
Let $\fg$ be  the graded Lie algebra  in Definition \ref{d.fh}.
Let $M$ be a simply connected complex manifold equipped with a regular distribution $\sD \subset TM$ and a fiber subbundle $\sS \subset \BP \sD$  such that  for each $y \in M$ \begin{itemize} \item[(1)] there exists an isomorphism of graded Lie algebras $\varphi_0: {\rm Symb}_y \sD \to \fg_+$  and \item[(2)]  the isomorphism $\varphi_0$ sends the projective subvariety $\sS_y \subset \BP \sD_y$ isomorphically to the Segre variety $\BP \bW \times \BP \fl_1 \subset  \BP \fg_1$.  \end{itemize}
Let  $\widetilde{\sH}_y$ for each  $y \in M$ be  the set of graded Lie algebra isomorphisms from
 $\fg_{+}$ to ${\rm Symb}_y \sD$  that sends  $\BP \bW \times \BP \fl_1$ to $\sS_y$.  Then $\widetilde{\sH} = \cup_{y \in M} \widetilde{\sH}_y$ is a $\widetilde{\bG_0}$-principal bundle on $M$ where $\widetilde{\bG}_0$ is as in Proposition \ref{p.auttimes}. Since $M$ is simply connected, we can  choose a connected component of $\widetilde{\sH}$ to have a $\bG_0$-principal bundle $\sH$ on $M$.
 For a $\bG_0$-module $V$, denote by $\sE_M(V)$ the vector bundle associated with the  $\bG_0$-principal bundle $\sH \to M$. \end{setup}

The model of Setup \ref{setup} is the following example.

\begin{example}\label{e.model}
In Definition \ref{d.prolong}, let  $\bG$ be a connected Lie group of adjoint type with Lie algebra $\fg$ and $\bG^0 \subset \bG$ be the connected closed subgroup with Lie algebra $\oplus_{k \leq 0} \fg_k$. The homogeneous space $M^{\rm model}:=\bG/\bG^0$ has the $\bG^0$-principal bundle $\bG \to \bG/\bG^0$ and the action of $\bG^0$ on $\bg_1$ and $\BP \bW \times \BP\fl_1 \subset \BP \bg_1$ induce a distribution $\sD^{\rm model} \subset TM^{\rm model}$  and a fiber subbundle $\sS^{\rm model} \subset \BP \sD^{\rm model}$   on $M^{\rm model}$.  \end{example}

We have the following consequence of Theorem \ref{t.HL}.

\begin{theorem}\label{t.HL2}
In the setting of Setup \ref{setup}, assume $H^0(M, \sE_M(Q^k(\fg))) =0$  for any negative integer $k$. Then there exists an open immersion $\phi: M \to M^{\rm model}$ that sends $\sD$ (resp. $\sS$) isomorphically to $ \sD^{\rm model}|_{\phi(M)}$ (resp. $\sS^{\rm model}|_{\phi(M)}$)  in Example \ref{e.model}.
\end{theorem}

In the rest of the section, we show that the condition $H^0(M, \sE_M(Q^k(\fg))) =0$ holds under the following additional assumption.

 \begin{assumption}\label{assumption}
In Setup \ref{setup},
let $\{\exp(tv_0) \mid t \in \C\}$ be the 1-parameter subgroup of $\bG_0$ corresponding to $v_0 \in \bV_0 \subset \fg_0$ and consider the adjoint action
$$\fg_+ \stackrel{\exp(tv_0)}{\longrightarrow} \fg_+$$  of each element $\exp(t v_0) \in \bG_0$. Define $\varphi_t := \exp(t v_0) \cdot \varphi_0$  as the composition $$ {\rm Symb}_y \sD \stackrel{\varphi_0}{\longrightarrow} \fg_+ \stackrel{\exp(t v_o)}{\longrightarrow} \fg_+,$$ where $\varphi_0$ is as in Setup \ref{setup}.
We assume that for a general point $y \in M$,  there exist  families $$ \{ \psi_t: Z \to M \mid t \in \Delta\} \mbox{ and } \{\zeta_t: \psi_t^* \sE_M(\fg_+) \to \sE_Z(\fg_+) \mid t \in \Delta\}$$ parametrized by $t \in \Delta$ for a disc $\Delta \subset \C$  centered at $0 \in \C,$ where each $\psi_t$ is a  holomorphic immersion and  each $\zeta_t$ is an  isomorphism of  vector bundles on $Z$  satisfying   the following for each $t \in \Delta.$
\begin{itemize} \item[(1)]  $\psi_t(z= [v_0]) = y$;  \item[(2)] ${\rm d} \psi_t (T_z Z) \subset \sD_y$; \item[(3)] $\varphi_t ({\rm d} \psi_t (T_z Z)) =  \fl_1 \subset \fg_1$; and
\item[(4)] the composition of $\varphi_0$ and the restriction of $\zeta_t$ to the fibers at $z= [v_0] \in Z,$
    $$ {\rm Symb}_y \sD \stackrel{\varphi_0}{\longrightarrow} \sE_M(\fg_+)_y  \stackrel{\zeta_t}{\longrightarrow}  \sE_Z (\fg_+)_z =  \fg_+, $$
    coincides with $\varphi_t.$
\end{itemize} \end{assumption}

\begin{lemma}\label{l.repara}
Consider the following diagram given by Assumption \ref{assumption},
 \[ \begin{array}{ccccc}
T_z Z & \stackrel{{\rm d} \psi_t}{\longrightarrow} & {\rm d} \psi_t (T_z Z) & \subset & \sD_y \\
\| & & \downarrow  & &\downarrow \\
\fl_1 & \stackrel{h_t}{\longrightarrow} & \fl_1& \subset & \fg_1, \end{array} \] for each $t \in \Delta$, where $h_t$ is an isomorphism of vector spaces and the two vertical arrows are isomorphisms given by $\varphi_t$. Then we may assume that $h_t = {\rm Id}_{\fl_1}$ for any $t \in \Delta$ by replacing $\psi_t$ by $\psi_t \circ \eta_t^{-1}$ for some family of elements $\{ \eta_t \in L_0  \mid t \in \Delta\}.$
\end{lemma}

\begin{proof} Let ${\rm Bs}(\III) \subset \BP T_z Z = \BP \fl_1$ be the base locus of the third fundamental form.
Since $h_t$ is the restriction of the isomorphism $\varphi_t$ of the graded Lie algebra $\fg_+,$ it sends ${\rm Bs}(\III) \subset \BP \fl_1$ to itself. This implies  $h_t \in L_0$. Thus we can find $\eta_t \in L_0$ with the required property.
\end{proof}

 \begin{lemma}\label{l.sH}
In Setup \ref{setup}, suppose that we have an immersion $\psi: Z \subset M$ and an isomorphism of vector bundles $$\zeta: \psi^* {\rm Symb} \sD \to \sE_Z(\fg_+)$$ which sends $\psi^* \sS$ to the fiber subbundle of
$\sE_Z (\fg_+)$ associated with $\BP \bW \times \BP \fl_1.$ Then $\zeta$
induces a morphism $\zeta_*: L \to \psi^* \sH$ of the $P$-principal bundle $L \to L/P = Z$ to the $\bG_0$-principal bundle $\psi^* \sH$ on $Z$ such that ${\rm Im}(\zeta_*)$  is  an $L_0$-principal subbundle  of  $\psi^*\sH$. Consequently, it induces an isomorphism $\psi^* \sE_M(V) \cong \sE_Z(V)$ for any $\bG_0$-module $V$. \end{lemma}

\begin{proof}
The morphism $ \zeta_*: L \to \psi^* \sH$ is defined as follows.
For $g \in L$, set $z = [gP] \in L/P =Z$ and define $\zeta_*(g) \in \sH_{\psi(z)}$ as the composition
$$\fg_+ = \sE_Z(\fg_+)_{[P]}  \stackrel{g}{\longrightarrow}  \sE_Z(\fg_+)_z \stackrel{\zeta^{-1}}{\longrightarrow}
 {\rm Symb}_{\psi(z)}\sD.$$
Since $\zeta_*(g) = \zeta_*(gp) $ exactly when  $p \in P$ belongs to the subgroup with Lie algebra $\fl_{-1}$, the image of $\zeta_*$ is an $L_0$-principal subbundle. \end{proof}

The following is a direct consequence of Lemma \ref{l.sH} applied to each $\zeta = \zeta_t$.

\begin{corollary}\label{c.sH}
In Setup \ref{setup} with Assumption \ref{assumption},
the isomorphism $\zeta_t$ induces an isomorphism $\zeta_t^V: \psi_t^* \sE_M(V) \cong \sE_Z(V)$ for any $\bG_0$-module $V$ for any $t \in \Delta$.
\end{corollary}

\begin{theorem}\label{t.vanishing}
In the setting of Setup \ref{setup} and Assumption \ref{assumption}, for any  negative integer $k$, we have
$$H^0(M, \sE_M(Q^{k}(\fg))) =0.$$
\end{theorem}

\begin{proof}
We may assume that we have made the choice of $\zeta_t$ in Assumption \ref{assumption} such that $h_t$ in Lemma \ref{l.repara} is ${\rm Id}_{\fl_1}$.
Choose any section $F \in H^0(M, \sE_M(Q^{k}(\fg)))$.
For a general point $y\in M$, we use the choice of $\psi_t$ as in Lemma \ref{l.repara} and define $F_t \in H^0(Z, \sE_Z(Q^{k}( \fg)))$ as the image of $$\psi_t^*F
\in H^0(Z, \psi_t^* \sE_M(Q^{k}(\fg))) \stackrel{\zeta_t^{Q^{k}(\fg)}}{\longrightarrow} H^0(Z, \sE_Z(Q^{k}(\fg))),$$ where $\zeta_t^{Q^k(\fg)}$ is as in Corollary \ref{c.sH}.
By Proposition \ref{p.R_k}, the value $F_{t, z}$ of $F_t$ at the base point $z \in Z$ is sent to an element $f_t \in R_k \subset \Hom(\wedge^2 \fg_+, \fg)_k$ such that $f_t=0$ if and only if $F_{t,z} =0$. Thus to prove the theorem, it suffices to check that $f_t =0$.

 Denote by $A$  the endomorphism of $\fg$ given by the Lie bracket $[v_o, \cdot].$  When regarding $\exp(sv_0)$ as an operator on $\fg$, we have the equalities $$ \exp(sv_0) = {\rm Id}_{\fg} + s A \mbox{ and } ({\rm Id}_{\fg} + sA)^{-1} = {\rm Id}_{\fg} -s A \  \mbox{ for } s \in \C.$$ By our assumption that $h_t = {\rm Id}_{\fl_1}$ in Lemma \ref{l.repara}, we have $f_{s+t} = \exp(sv_0) \cdot f_t$ for any $t\in \Delta$ and sufficiently small $s$, where the right hand side denotes the action of $\exp(sv_0)$ on $\Hom(
   \wedge^2 \fg_+, \fg)_k$. Thus  for any $u, u' \in \fg_+$, $t\in \Delta$ and sufficiently small  $s$,  we have $$f_{s+t}(u, u') = ({\rm Id}_{\fg} + sA) \cdot f_t(({\rm Id}_{\fg} - sA) u, ({\rm Id}_{\fg} -s A) u'),$$ which gives   \begin{equation}\label{e.st}  f_{s+t} (u, u')  =  
\begin{array}{l}  f_t(u, u')  +  s A \cdot f_t(u,u') \\ - s f_t(u, A \cdot u')   - s f_t( A \cdot u, u') \\ + s^2 f_t(A \cdot u, A \cdot u')  - s^2 A \cdot f_t(u, A \cdot u')  \\ - s^2 A\cdot f_t(A \cdot u, u')  + s^3 A \cdot f_t(A \cdot u, A \cdot u'). \end{array} 
\end{equation}
Using (\ref{e.st}), we  check $0= f_t \in R_k$ for $k = -1, -2, -3$ one by one in the following three lemmata. \end{proof}

\begin{lemma}\label{l.R-1}
For $f_t \in R_{-1} = \Hom(\bV_1 \wedge \bV_+, \bV)_{-1} \oplus \Hom (\fl_1 \wedge \bV_1, \fl_1)\oplus\Hom(\wedge\bV_1, \fl_1),$
we have \begin{itemize} \item[(1)] $f_t(\wedge^2 \bV_1) =0$;
\item[(2)] $f_t(\bV_1 \wedge \bV_2) = f_t(\bV_1 \wedge \bV_3) = 0$;
\item[(3)] $f_t(\fl_1 \wedge \bV_1) = 0$;
\end{itemize} which implies $f_t =0$. \end{lemma}

    \begin{proof}
    To check (1), choose $u, u' \in \fl_1 $ in (\ref{e.st}) to have
    $$0 = f_{s+t}(u, u')  = s F_{1, t}(u, u') + s^2 F_{2, t}(u,u') + s^3 A \cdot f_t(A\cdot u, A\cdot u'),$$
    where
    \begin{eqnarray*}
    && F_{1, t}(u, u'):=-f_t(u, A\cdot u')-f_t(A\cdot u, u'), \\
    && F_{2, t}(u, u'):=f_t(A\cdot u, A\cdot u')+A\cdot F_{1, t}(u, u').
    \end{eqnarray*}
    As this holds for all sufficiently small $s$, we have $$F_{1,t}(u,u') = 0 = F_{2,t}(u,u'),$$
    which implies $f_t(A \cdot u, A \cdot u') =0$.
This proves (1) because $A \cdot \fl_1 = [v_0, \fl_1] = \bV_1.$

    To check (2), choose $u \in \fl_1$ and $u' \in \bV_2 \oplus \bV_3$ in (\ref{e.st}) to have
    $$0 = f_{s+t}(u, u') = -s f_t(A \cdot u, u') + s^2 (\cdots) + s^3 (\cdots).$$
     As this holds for all sufficiently small $s$, we have $f_t(A \cdot u, u') =0$,
     which proves (2) because $A \cdot \fl_1 = [v_0, \fg_1] = \bV_1.$

     To check (3), choose $u \in \fl_1$ and $u' \in \bV_1$ in (\ref{e.st}) to have
     $$ \fl_1 \ni f_{s+t}(u, u') = f_t (u, u') +s A \cdot f_t(u, u') + s^2(\cdots) + s^3(\cdots),$$
     where we have used $f_t(\wedge^2\bV_1)=0$ from (1).
     As this holds for all sufficiently small $s$ and $f_t(u, u') \in \fl_1,$
     we have $$ \fl_1 \ni A \cdot f_t(u, u') = [v_0, f_t(u, u')]  \in [v_0, \fl_1] = \bV_1.$$
     Since $\fl_1 \cap \bV_1 =0$ as subspaces of $\fg_1$, we have $A\cdot f_t(u, u')=0,$ which proves $f_t(u, u')=0$.
    \end{proof}

     \begin{lemma}\label{l.R-2}
For $f_t \in R_{-2} =  \Hom(\wedge^2 \bV_+, \wg)_{-2} ,$
we have  $$f_t(\bV_1  \wedge (\bV_1 \oplus \bV_2)) = 0,$$ which implies $f_t =0$. \end{lemma}

     \begin{proof}
     Choose $u \in \fl_1$ and $u' \in \bV_1 \oplus \bV_2$ in (\ref{e.st}) to have
     $$0 = f_{s+t}(u, u') = -s f_t(A \cdot u, u') + s^2 (\cdots) + s^3 (\cdots) . $$
     As this holds for all sufficiently small $s$,
     we have $f_t( A \cdot u, u') =0$, which implies $f_t(\bV_1, \bV_1 \oplus \bV_2) =0$ because $\bV_1 = A \cdot \fl_1.$ \end{proof}

     \begin{lemma}\label{l.R-3}
     For $f_t \in R_{-3} =  \Hom(\wedge^2 \bV_+, \wg)_{-3},$
we have \begin{itemize} \item[(1)] $f_t(\bV_1 \wedge \bV_+) =0$;
\item[(2)] $f_t(\wedge^2 \bV_2) =0,$ \end{itemize} which implies $f_t =0$.  \end{lemma}

     \begin{proof}
     To check (1), choose $u \in \fl_1, u' \in \bV_+$ in (\ref{e.st}) to have $$0 = f_{s+t}(u, u') = -s f_t(A \cdot u, u') + s^2( \cdots) + s^3 (\cdots).$$
     As this holds for all sufficiently small $s$, we have $f_t(A \cdot u, u') =0$, which implies (1) because $\bV_1 = A \cdot \fl_1$.

     To check (2), choose $u, u' \in \bV_2$ in (\ref{e.st}) and use $A \cdot u = A \cdot u'=0$ to have $$f_{s+t}(u, u') - f_t(u, u') = s A \cdot f_t(u, u').$$ Since the left hand side is in $\fl_1$ and the right hand side is in $\bV_1$, we obtain $A\cdot f_t(u,u') =0,$ which implies $f_t(u,u') =0$.
\end{proof}

\section{Minimal rational curves whose VMRT's are subadjoint varieties}\label{s.vmrt}

Let us recall the notion of the variety of minimal rational
tangents.

\begin{definition}\label{d.VMRT}
Let $X$ be a uniruled projective manifold. For an irreducible
component  $\sK$ of the space of rational curves on $X$, denote by
$\rho: {\rm Univ}_{\sK} \to \sK$ and $\mu: {\rm Univ}_{\sK} \to X$ the universal family
morphisms. The component $\sK$ is called a {\em family of minimal rational curves}, if for a general $x \in X$, the fiber $\mu^{-1}(x)$ is
non-empty and projective (or compact). In this case, replacing ${\rm Univ}_{\sK}$ by its normalization if
necessary,  we can assume that $\rho$ is a $\BP^1$-bundle and
$\sK_x := \mu^{-1}(x)$ is a smooth projective variety for a general
$x \in X$ (e.g. Theorem 1.3 of \cite{Hw01}).    The {\em tangent map} $\tau: {\rm Univ}_{\sK} \dasharrow \BP TX$ is the rational map
associating  a smooth point of a rational curve to its tangent
direction. For a family of minimal rational curves, the tangent map
induces a morphism $\tau_x: \sK_x \to \BP T_x X$ for a general point $x \in X$, which is a normalization of its image $\sC_x \subset \BP T_x X$ (by Theorem 3.4 of \cite{Ke} and Theorem 1 of \cite{HM04}). This image $\sC_x$ is the {\em variety of minimal
rational tangents} (VMRT) at $x$. Denote by $\sC \subset \BP
TX$ the proper image of $\tau$.
\end{definition}

Throughout this section, we work in the following setting.

\begin{setup}\label{setupX}
Fix a subadjoint variety $Z \subset \BP \bV$ excluding the case (0) of Proposition \ref{p.list}.
Let $X$ be a uniruled projective manifold with a family $\sK$ of minimal rational curves and a nonempty open subset $X^o \subset X$ such that the VMRT $\sC_x \subset \BP T_x X$ at every $x \in X^o$ is projectively isomorphic to $Z \subset \BP \bV.$ Thus $\sC^o := \cup_{x \in X^o} \sC_x$ is a $Z$-isotrivial cone structure on $X^o$. Since  the tangent morphism $\tau_x: \sK_x \to \sC_x$ is the normalization of $\sC_x$ for a general $x \in X$, we may assume that $\tau_x$ is biregular  for every $x \in X^o.$  In particular, the universal morphism $\rho: {\rm Univ}_{\sK} \to \sK$ induces a  holomorphic submersion $\varrho: \sC^o \to \sK^o$ to an open subset $\sK^o$ in the smooth locus of $\sK$ whose fibers are transversal to fibers of the projection $\pi: \sC^o \to X^o$ (see page 58 of \cite{HM04} for more details).  For  $\kappa \in \sK^o$, denote by  $f_{\kappa}: \BP^1 \to X$ the restriction of the universal morphism $\mu:{\rm Univ}_{\sK} \to X$ to $\rho^{-1}(\kappa)$ composed with an identification $\BP^1 \cong \rho^{-1}(\kappa)$.  Then $f_{\kappa}$ gives a normalization of its image and satisfies
$$f_{\kappa}^* TX \cong \sO(2) \oplus \sO(1)^{\oplus m} \oplus \sO^{\oplus (n-1-m)}$$
for $m = \dim \sC_x = \dim Z$ and $n= \dim X$ (e.g. by Proposition 1.4 and its proof in \cite{Hw01}). From Proposition \ref{p.FFsubadjoint}, we know $n = 2m+2$. \end{setup}

\begin{notation}\label{n.cone}
In the setting of Setup \ref{setupX}, we have vector subbundles $\sV \subset  \sT^0 \subset \sT^1 \subset \sT^2 \subset T \sC^o$ such that their fibers at $z \in \sC_x \subset \sC^o$ with $x \in X^o$ are
\begin{eqnarray*} \sV_z & = & ({\rm d}_z \pi)^{-1} (0) \\
\sT^0_z &=& ({\rm d}_z \pi)^{-1} (\widehat{z}) \\
\sT^1_z & = & ({\rm d}_z \pi)^{-1} (\widehat{T}_z \sC_x)\\
\sT^2_z & = & ({\rm d}_z \pi)^{-1}(T^{(2)}_z \sC_x) \end{eqnarray*} where ${\rm d}_z \pi: T_z \sC^o \to T_x X^o$ is the differential of the projection $\pi: \sC^o \to X^o$ at $z$.
\end{notation}

\begin{proposition}\label{p.D}
In Setup \ref{setupX}, recall (see p. 58 of \cite{HM04}) that the tangent space $T_{\kappa} \sK^o$ can be naturally identified with $$H^0(\BP^1, f_{\kappa}^* TX)/H^0(\BP^1, T \BP^1)  \cong H^0(\BP^1, \sO(1)^{\oplus m} \oplus \sO^{\oplus (m+1)}),$$ where $m= \dim Z$.  Define a vector subbundle $D \subset T \sK^o$
whose fiber at $\kappa$ corresponds to the subspace $H^0(\BP^1, \sO(1)^{\oplus m})$ of the right hand side in the above isomorphism.
Then \begin{itemize} \item[(1)] $\sT^1_z = ({\rm d}_z \varrho)^{-1} D_{\varrho(z)}$ for each $z \in \sC^o$; \item[(2)] the symbol algebra ${\rm Symb}_{\kappa} D $ is isomorphic to $\fg_+$; and  \item[(3)] there exists a subbundle $J \subset D$ of rank $m$ corresponding to $\bV_1 \subset \fg_1$ under the isomorphism in (2) such that $\sT_z^1 = \sV_z \oplus ({\rm d}_z \varrho)^{-1} J_{\varrho(z)}.$ \end{itemize} \end{proposition}

\begin{proof}
(1) is just Proposition 8 in \cite{HM04}.  By (1), the symbol algebra of $D$ at $\kappa$ is isomorphic to the symbol algebra of $\sT^1$ at a point of $\varrho^{-1}(\kappa)$, modulo $\sF := {\rm Ker}({\rm d} \varrho).$
We have a natural isomorphism \begin{equation}\label{e.otimes} \sV \otimes \sF  \ \cong \ \sT^1/\sT^0 \end{equation} by Proposition 1 of \cite{HM04}.  At a point $z \in \varrho^{-1}(\kappa)$, write \begin{eqnarray*} \bg_1 &:=& \sT^1_{z}, \\ \bg_2 &:=& \sT^2_{z}/\sT^1_{z}, \\ \bg_3 &:=& \sT^3_{z}/\sT^2_{z}. \end{eqnarray*}
We have a natural inclusion $$ \sV_{z} \oplus \sF_{z} = \sT_{z}^0  \ \subset \sT^1_{z}$$ and  a natural isomorphism $\xi: \sV_{z}  \to \sT^1_{z}/\sT^0_{z}$ determined by a choice of a nonzero vector in $\sF_{z}$ in (\ref{e.otimes}).
The symbol algebra of $\sT^1$ is $\bg_+ = \bg_1 \oplus \bg_2 \oplus \bg_3$ with the Lie brackets satisfying (by Proposition 3.16 of \cite{HL2})
$$[\sF_{z}, \sT_{z}^0] =0, \   [\sV_{z}, \sV_{z}] = 0$$ and for any $u, v, w \in \sV_{z}$ $$[u, \xi(v)] = \II(u, v), \ [w, \II(u, v)] = \III(u,v,w) $$ where $\II$ and $\III$ are the second and the third fundamental forms of $\sC_x \subset \BP T_x X$ at the point $z$.
By Proposition \ref{p.FFsubadjoint} (3),  the bilinear map $[\sV_{\alpha}, \bg_2] \to \bg_3$ determines a perfect pairing between $\sV_{z}$ and $\bg_2$. It follows (see Lemma \ref{l.pairing} below) that there exists a subspace $W \subset \sT^1_{z}$ complementary to $\sT^0_{z}$ such that $[W, \bg_2] =0$.
The isomorphisms $$\sV_{z} \stackrel{\xi}{\cong} \sT^1_{z}/\sT^0_{z}\cong W$$  induce an isomorphism
between $\sV_{\alpha} \cong \fl_1$ and $ W \cong \bV_1$.  Together with $ \bg_2 \cong \bV_2$ and  $\bg_3 \cong \bV_3$, it gives the isomorphism between ${\rm Symb}_{\kappa} D$ and $\fg_+$, which proves (2). The image $J_{\varrho(z)}$ of $W$ in $T_{\varrho(z)} \sK^o$ depends only on $\kappa = \varrho(z)$, because it is determined by ${\rm Symb}_{\kappa} D$. Thus we have the vector bundle $J$ of (3).
\end{proof}

\begin{lemma}\label{l.pairing}
Let $\gamma: A \otimes B \to \C$ be a bilinear map. Assume that there exists a subspace $B' \subset A$ such that the induced bilinear map $B' \otimes B \to \C$ is a perfect pairing. Then there exists a subspace $W \subset A$ complementary to $B'$ such that $\gamma(W, B) = 0$. \end{lemma}

\begin{proof}
The homomorphism $j: B \to A^{\vee}$ induced by $\gamma$ is injective because $B' \otimes B \to \C$ is a perfect pairing. Let $W \subset A$ be the subspace annihilated by $j(B)$. Then $W$ satisfies the required properties.
\end{proof}

\begin{proposition}\label{p.sS}
In Proposition \ref{p.D}, for each $\kappa\in \sK^{o}$  and $f_{\kappa}: \BP^1 \to X$ with  a point $t \in \BP^1$ satisfying $f_{\kappa}(t) =: x_t \in X^o$, there exists a natural isomorphism \begin{equation}\label{e.tensor} D_{\kappa} \cong H^0(\BP^1, \sO(t)) \otimes H^0(\BP^1, (f_{\kappa}^* TX/T\BP^1) \otimes {\bf m}_t), \end{equation} where $H^0(\BP^1, \sO(t))$ is the space of rational functions on $\BP^1$ with poles at $t$ and ${\bf m}_t$ denotes the sheaf of ideals of local holomorphic functions on $\BP^1$ vanishing at $t$. Writing $z_t = \varrho^{-1}(\kappa) \cap \sC_{x_t}$, this induces a natural identification  $$H^0(\BP^1, (f^* TX/T\BP^1) \otimes {\bf m}_t) = {\rm d}_{z_t} \varrho (T_{z_t} \sC_{x_t}) \ \subset \ T_{\kappa} \sK^o.$$  Moreover, if we define $\sS_{\kappa} \subset \BP D_{\kappa}$ as the projectivization of the set of decomposable tensors in the tensor decomposition (\ref{e.tensor}), then $\sS_{\kappa}$ is  projectively isomorphic to $\BP \bW \times \BP \fl_1 \subset \BP (\bW \otimes \fl_1)$ and depends only on $\kappa$,  independent of $t$. \end{proposition}

\begin{proof}
All are straightforward. Let us just check the last statement that $\sS_{\kappa}$ is independent of $t$.
If we choose another point  $s \in \BP^1$  different from $t$ satisfying $ f_{\kappa}(s) \in X^o$, the tensor decomposition is
$$D_{\kappa} \cong H^0(\BP^1, \sO(s)) \otimes H^0(\BP^1, (f^* TX/T\BP^1) \otimes {\bf m}_{s}).$$
Fix a rational function $h$ on $\BP^1$ with the zero divisor $t$ and the pole divisor $s$.
Then $H^0(\BP^1, \sO(s)) = h \cdot H^0(\BP^1, \sO(t))$ and $$H^0(\BP^1, (f^* TX/T\BP^1) \otimes {\bf m}_t) = h \cdot H^0(\BP^1, (f^* TX/T\BP^1) \otimes {\bf m}_{s}).$$
Thus the set of decomposable tensors do not depend on the choice of $t$ and the submanifold  $\sS_{\kappa} \subset \BP D_{\kappa}$ is determined by $\kappa \in \sK$.  \end{proof}

\begin{proposition}\label{p.zeta}
In Proposition \ref{p.sS}, fix a biholomorphic identification of a complex disc $\Delta \subset \C$ with an open subset in $\BP^1 \cap f_{\kappa}^{-1}(X^o)$ and choose a holomorphic family of biregular morphisms $\{\Psi_t: Z \cong \sC_x\}$.  Set $\psi_t = \varrho \circ \Psi_t : Z \to \sK^o,$ which is an immersion. Then, after shrinking $\Delta$ if necessary,  there exists a holomorphic family of  isomorphisms of vector bundles $$\{ \zeta_t: \psi_t^* {\rm Symb} D \to \sE_Z(\fg_+) \mid t \in \Delta\},$$ such that each $\zeta_t$ sends $\psi_t^*\sS$ to the fiber subbundle of $\BP \sE_Z(\fg_+)$ associated with $\BP \bW \times \BP \fl_1.$ \end{proposition}

\begin{proof}
For any $x \in X^o$, the line bundle $\sF := {\rm Ker}({\rm d} \varrho)$ restricted to $\sC_x$ is the tautological line bundle of $\sC_x \subset \BP T_x X$. Thus \begin{equation}\label{e.1} \Psi_t^* \sF \cong \sE_Z(\bV_0). \end{equation} From $\sV|_{{\sC}_x} = T\sC_x$, we have \begin{equation}\label{e.2} \Psi_t^* \sV \cong \sE_Z(\fl_1). \end{equation}
From (\ref{e.otimes}), (\ref{e.1}) and (\ref{e.2}),  we have \begin{equation}\label{e.3} \Psi_t^* (\sT^1/\sT^0) \cong \sE_Z (\bV_0 \otimes \fl_1) \cong \sE_Z(\bV_1). \end{equation}
From $(\sT^2/\sT^1)|_{\sC_x} \cong \sF \otimes N^{(2)}_{\sC_x}, N^{(2)}_Z \cong \sE_Z(N^{(2)}_{Z,z})$ and (\ref{e.1}), we have \begin{equation}\label{e.4} \Psi_t^* (\sT^2/\sT^1) \cong \sE_Z (\bV_0 \otimes N^{(2)}_{Z,z}) \cong \sE_Z(\bV_2). \end{equation} Similarly, we have \begin{equation}\label{e.5}  \Psi_t^* (\sT\sC^o/\sT^2) \cong \sE_Z (\bV_0 \otimes N^{(3)}_{Z,z}) \cong \sE_Z(\bV_3). \end{equation}
By Proposition \ref{p.D} (1), we have $$\varrho^*{\rm Symb} D \cong \sT^1/\sF \oplus \sT^2/\sT^1 \oplus \sT\sC^o/\sT^2$$ and $\sT^1/\sF \cong \sV \oplus (\sV \otimes \sF).$  Thus (\ref{e.3}), (\ref{e.4}) and (\ref{e.5}) show that   the vector bundles  $\psi_t^* {\rm Symb} D$ and $\sE_Z(\fg_+)$ on $Z$ are isomorphic.
By shrinking $\Delta$ if necessary, we can choose the family of isomorphisms $\zeta_t$.  The isomorphisms are compatible with the decompositions $\sT^1/\sF \cong \sV \oplus (\sV \otimes \sF)$ and $\fg_1 = \fl_1 \oplus (\fl_1 \otimes \bV_0).$ Thus it sends $\psi_t^* \sS$ to the fiber subbundle associated with $\BP \bW \times \BP \fl_1$.\end{proof}

We are ready to prove the following precise version of Theorem \ref{t.main}.

\begin{theorem}\label{t.main2}
In Setup \ref{setupX}, the $Z$-isotrivial cone structure $\sC^o \subset \BP TX^o$ is locally flat. \end{theorem}

To prove the theorem, it is convenient to use the following terminology.

\begin{definition}\label{d.preserve}
Recall that a vector field $\vec{v}$ on a complex manifold $Y$ can be lifted to a vector field on $\BP TY$: the local 1-parameter family of germs of biholomorphic maps of $Y$ generated by $\vec{v}$ induces a local 1-parameter family of germs of biholomorphic maps of $\BP TY$ and its derivative gives the lifted vector field on $\BP TY$.
Given a cone structure $\sC \subset \BP TY$ on a complex manifold $Y$, we say that a holomorphic vector field $\vec{v}$ on $Y$ {\em preserves } the cone structure, if the  lift of $\vec{v}$ to $\BP TY$ is tangent to $\sC$. \end{definition}

We have the following criterion of local flatness of a cone structure by Propositions 5.10 and 5.12 of \cite{FH12}.

\begin{lemma}\label{l.flat}
Let $\sC \subset \BP TY$ be a cone structure on a complex manifold $Y$ such  that for  a general $y \in Y$, the prolongation of $\aut(\widehat{\sC}_y) \subset \Hom(T_y Y, T_y Y),$ i.e. the vector space
$$\{ \varphi \in \Hom(\Sym^2 T_yY, T_y Y) \mid \varphi(v, \cdot) \in \aut(\widehat{\sC}_y) \mbox{ for any } v \in T_y Y \},$$ is zero.
Suppose there exists a Lie algebra $\fh$ of holomorphic vector fields on $Y$ preserving the cone structure with $\dim \fh = \dim Y + \dim \aut(\widehat{\sC}_y)$ for a general point $y \in Y$. Then $\sC$ is locally flat. \end{lemma}

\begin{proof}[Proof of Theorem \ref{t.main2}]
Fix a base point $x \in X^o$. Recall that  the fibers of $\pi: \sC^o \to X^o$ are simply connected because $Z$ is simply connected. Thus we can find a neighborhood $O \subset X^o$ of $x$ such that $\pi^{-1}(O)$ is simply connected.   The fibers of $\varrho: \sC^o \to \sK^o$ form a foliation on $\pi^{-1}(O)$ which is transversal to the fibers of $\pi: \sC^o \to X^o$. By shrinking $O$ if necessary, we can assume that the space $M$ of leaves of this foliation is a complex manifold such that \begin{itemize}
\item[(1)] the restriction $\varrho|_{\pi^{-1}(O)}$ defines a holomorphic submersion with connected fibers $\rho_O: \pi^{-1}(O) \to M$ onto the  complex manifold $M$; \item[(2)] for any $y \in O$, the holomorphic map $\rho_O$  sends $\sC_y$  biholomorphically to a submanifold in $M$; and \item[(3)] there exists an open immersion $j: M \to \sK^o$ such that $\varrho|_{\pi^{-1}(O)} = j \circ \rho_O$. \end{itemize}  Note that $M$ is simply connected because $\pi^{-1}(O)$ is simply connected.

     When $n= \dim X$, we have the exact sequence for each $y \in O$ (as in Proposition 4 of \cite{HMaut}) $$0 \to \sF|_{\sC_y} \to \sO^{\oplus n}_{\sC_y} \to \rho_O^* N_{\rho_O(\sC_y)} \to 0$$ where $\sF = {\rm Ker}({\rm d} \rho_O)$ and  $N_{\rho_O(\sC_y)}$ is the normal bundle of the submanifold $\rho_O(\sC_y) \subset M$.
Since the line bundle $\sF|_{\sC_y}$ is isomorphic to the restriction of $\sO(-1)_{\BP \bV}$ to $Z$,
Kodaira vanishing theorem gives $\dim H^0(\rho_O(\sC_y), N_{\rho_O(\sC_y)}) = n$. This means that all small deformations of the submanifold $\rho_O(\sC_x)$ in $M$ come from images of $\sC_y, y \in O$. Consequently, the pair of holomorphic maps \begin{equation}\label{e.and} \rho_O:\pi^{-1}(O) \to M \mbox{ and } \pi|_{\pi^{-1}(O)}: \pi^{-1}(O) \to O \end{equation} can be viewed as the universal family of small deformations of $\rho_O(\sC_x)$ in $M$. More precisely, there is an open embedding of $O$ into the Douady space of $M$ (in the sense of Section VIII.1 of \cite{GP}) such that $\pi|_{\pi^{-1}(O)}$ is the pullback of the universal family on the Douady space.

The immersion $j: M \to \sK^o$  induces a distribution $\sD$ on $M$ by pulling back $D \subset T\sK^o$ and a fiber subbundle $\sS \subset \BP \sD$ given
by $\sS_{\kappa}, \kappa \in j(M) \subset \sK^o,$ from Proposition \ref{p.sS}.  They satisfy the conditions in Setup \ref{setup}. By Proposition \ref{p.zeta}, they satisfy Assumption \ref{assumption}. Then by Theorem \ref{t.vanishing},  the condition in Theorem \ref{t.HL2} is satisfied. Thus we have an open immersion $\phi: M \to \bG/\bG_0$.

The holomorphic vector fields on $\bG/\bG_0$ corresponding to $\fg$ can be pulled back to holomorphic vector fields on $M$ by $\phi$. Since (\ref{e.and}) comes from the universal family on the Douady space ${\rm Douady}(M)$ of $M$, these vector fields on $M$ induce holomorphic vector fields on  $\pi^{-1}(O)$ and $O$ by the universal property of the Douady space. In fact, for a vector field $v$ on $M$ and a compact analytic subvariety $ C \subset M$, the vector field $v$  generates a local 1-parameter family of open embeddings of a neighborhood of $C \subset M$ into $M$, hence a local 1-parameter family of open embeddings of a neighborhood of  the corresponding point   $[C] \in {\rm Douday}(M)$ into ${\rm Doudady}(M)$. The derivative of this 1-parameter family of embeddings gives the  vector field on the Douady space induced by $v$.
Because they are induced by the universal family, these vector fields on $O$ preserve  the cone structure $\sC|_O$ in the sense of Definition \ref{d.preserve} (This argument is an infinitesimal version of Theorem 1 of \cite{HMaut}). The Lie algebra of these vector fields on  $O$ is
 isomorphic to $\fg$.  Note that $\aut(\widehat{Z})$ has no nonzero prolongation (as it does not belong to  the list in Table 5 of \cite{MS}) and $\dim \fg = \dim O + \dim \aut(\widehat{Z})$ from Definition \ref{d.fh}.  Thus  the cone structure $\sC$ must be locally flat by Lemma \ref{l.flat}. \end{proof}

\bigskip
{\bf Acknowledgment}
We would like to thank Lorenz Schwachh\"ofer  for helpful discussion on Theorem \ref{t.MS} and the referee for helpful suggestions to improve the presentation of the paper.

\bigskip
Jun-Muk Hwang (jmhwang@ibs.re.kr)
Center for Complex Geometry,
Institute for Basic Science (IBS),
Daejeon 34126, Republic of Korea

\medskip
Qifeng Li (qifengli@amss.ac.cn)
School of Mathematics, 
Shandong University,
Jinan, China


\begin{thebibliography}{KSWZY}
\bibitem{Ak} Akhiezer, D.: {\rm Lie group actions in complex analysis}. Vieweg Verlag. Braunschweig/Wiesbaden, 1995
    \bibitem{BF} Brion, M. and Fu, B.:
    Minimal rational curves on wonderful group compactifications. J. Ec. polytech. Math. {\bf 2} (2015) 153-170
    \bibitem{Br} Bryant, R.: Metrics with exceptional holonomy. Annals Math. {\bf 126} (1987) 525-576

\bibitem{Bu} Buczynski, J.:  Legendrian subvarieties of projective space. Geom. Dedicata {\bf 118} (2006) 87-103

\bibitem{FH12} Fu, B. and Hwang, J.-M.: Classification of non-degenerate projective varieties with non-zero prolongation and application to target rigidity.  Invent. math. {\bf 189} (2012) 457-513
  \bibitem{FH18} Fu, B. and Hwang, J.-M.: Isotrivial VMRT-structures of complete intersection type.  Asian J. Math. {\bf 22}, {\it Special issue for Ngaiming Mok's 60th birthday}, (2018) 333-354


\bibitem{GP}  Grauert, H., Peternell, T. and Remmert, R.: {\em Several complex variables VII}. Encycl. Math. Sci. {\bf 74},  Springer-Verlag, Berlin-Heidelberg, 1994
\bibitem{Hum} Humphreys, J.: {\em Introduction to Lie algebras and representation theory}.
Graduate Texts in Mathematics. 9, Springer-Verlag, New York, 1972
\bibitem{Hw01} Hwang, J.-M.: Geometry of minimal rational curves
on Fano manifolds. {\em School on Vanishing Theorems and Effective
Results in Algebraic Geometry (Trieste, 2000)}, 335--393, ICTP
Lect. Notes, 6, Abdus Salam Int. Cent. Theoret. Phys., Trieste,
2001
\bibitem{Hw10} Hwang, J.-M.: Equivalence problem for minimal rational curves with isotrivial
varieties of minimal rational tangents.  Ann. scient. Ec. Norm.
Sup. {\bf 43} (2010) 607-620
\bibitem{Hw12}  Hwang, J.-M.:
Geometry of varieties of minimal rational tangents. in {\em
Current Developments in Algebraic Geometry}, Mathematical Sciences
Research Institute Publications {\bf 59} (2012) 197-226, Cambridge
University Press
\bibitem{Hw13}  Hwang, J.-M.:
Varieties of minimal rational tangents of codimension 1.
Ann. scient. Ec. Norm. Sup. {\bf 46} (2013) 629-649
\bibitem{Hw21}  Hwang, J.-M.: Legendrian cone structures and contact prolongations. in {\em Geometry, Lie Theory and Applications- The Abel symposium 2019},  Abel Symposia {\bf 16} (2022) 131-145

\bibitem{HL} Hwang, J.-M. and Li, Q.: Characterizing symplectic Grassmannians by varieties of minimal rational tangents. J. Differential Geom. {\bf 119} (2021) 309-381.

\bibitem{HL2} Hwang, J.-M. and Li, Q.: Unbendable rational curves of Goursat type and Cartan type. J. Math. Pures Appl. {\bf 155} (2021) 1-31.

\bibitem{HM99} Hwang, J.-M. and Mok, N.: Varieties of minimal
rational tangents on uniruled projective manifolds. {\em Several
complex variables (Berkeley, CA, 1995--1996)}, 351-389, Math.
Sci. Res. Inst. Publ., 37, Cambridge Univ. Press, Cambridge, 1999
\bibitem{HMaut}  Hwang, J.-M. and Mok, N.:
Automorphism groups of spaces of minimal rational curves on Fano manifolds of Picard number 1. J. Algebraic Geom. {\bf 13} (2004)  663-673

\bibitem{HM04}
Hwang, J.-M. and Mok, N.: Birationality of the tangent map
for minimal rational curves. Asian J. Math. {\bf 8} (2004) 51-63

\bibitem{HN} Hwang, J.-M. and Neusser, K.: Cone structures and parabolic geometries.
to appear in Math. Annalen,  arXiv:2010.14958


\bibitem{Ke}  Kebekus, S.: Families of singular rational curves. J. Alg. Geom. {\bf 11} (2002) 245-256

\bibitem{LM03} Landsberg, J. M. and Manivel, L.: On the projective geometry of rational homogeneous varieties. Comm. Math. Helv. {\bf 78} (2003) 65-100

    \bibitem{LM07} Landsberg, J. M. and Manivel, L.: Legendrian varieties. Asian J. Math. {\bf 11} (2007) 341-360


\bibitem{MS} Merkulov, S. and Schwachh\"ofer, L. : Classification of irreducible holonomies of torsion-free affine connections. Annals Math. {\bf 150} (1999) 77-149

\bibitem{Mk} Mok, N.: Recognizing certain rational homogeneous
manifolds of Picard number 1 from their varieties of minimal
rational
 tangents.  Third International Congress of Chinese Mathematicians.
 Part 1, 2, 41-61, AMS/IP Stud. Adv. Math., 42, pt.1, 2,
  Amer. Math. Soc., Providence, RI, 2008


\bibitem{Ya} Yamaguchi, K.: Differential systems associated with simple graded Lie algebras.
Adv. Study Pure Math. {\bf 22}  (1993) 413-494

\bibitem{Za} Zak, F. L.: Tangents and secants of algebraic varieties.  Transl. Math. Monographs {\bf 127}, Amer. Math. Soc., Providence, RI, 1993
\end{thebibliography}
\end{document}